\documentclass[12 pt]{amsart}
\usepackage{amsthm, amsfonts, amssymb, color}
\usepackage{mathrsfs}
\usepackage{amsmath}
\usepackage{amstext, amsxtra}
\usepackage{txfonts}
\usepackage{soul} 
\usepackage[colorlinks, linkcolor=black, citecolor=blue, hypertexnames=false]{hyperref}
\usepackage{comment}

\allowdisplaybreaks
\usepackage{geometry}
\newtheorem{theorem}{Theorem}[section]
\newtheorem{proposition}[theorem]{Proposition}
\newtheorem{corollary}[theorem]{Corollary}
\newtheorem{lemma}[theorem]{Lemma}

\newtheorem{assumption}{Assumption}[section]
\newtheorem{remark}[theorem]{Remark}

\newcommand{\s}{ L}
\newcommand{\eq}{\begin{equation}}
\newcommand{\eeq}{\end{equation}}
\newcommand{\eqa}{\begin{eqnarray}}
\newcommand{\eeqa}{\end{eqnarray}}
\newcommand{\tphi}{\tilde{\phi}}

\newcommand{\jBra}[1]{\left\langle #1 \right\rangle}

\numberwithin{equation}{section}

\title   [] {On the Existence of Self-Similar solutions for some Nonlinear Schr\"odinger equations}

\author{   Avy Soffer }

\address{Department of Mathematics\\
Rutgers University\\
110 Frelinghuysen Rd.\\
Piscataway, NJ, 08854, USA}

\email{soffer@math.rutgers.edu}

\author{Xiaoxu Wu}
\address{Mathematical Sciences Institute\\ Australia National University\\ 
Canberra 2601, Australia.}
 \email{xiaoxu.wu@anu.edu.au}

\thanks{2010 \textit{ Mathematics Subject Classification.}   35Q55.  }
\thanks{
A.Soffer is supported in part by Simons Foundation Grant number 851844
}

\begin{document}

\begin{abstract}
We construct solutions of Schr\"odinger equations which are asymptotically self-similar solutions as time goes to infinity.
Also included are situations with two bubbles.
These solutions are global, with non-zero $L^2$ norms, and are stable.
As such they are not of the standard asymptotic decomposition of linear waves and localized waves.
Such weakly localized solutions were expected in view of previous works \cite{Liu-Sof1,Liu-Sof2} on the large time behavior of general dispersive equations.
It is shown that one can associate a \emph{scattering channel} to such solutions, with the dilation operator as the asymptotic ``Hamiltonian''.
\end{abstract}
\maketitle
\section{Introduction}
The analysis of dispersive wave equations and systems is of critical importance in the study of evolution equations in physics and geometry. It is well known that the asymptotic solutions of such equations, if they exist, show a dizzying zoo of possible solutions. Besides the ``free wave'', which corresponds to a solution of the equation without interaction terms, a multitude of other solutions may appear. Such solutions are localized around possibly moving centers of mass. They include nonlinear bound states, solitons, breathers, hedgehogs, vortices, etc. The analysis of such equations is usually done on a case by case basis, due to this complexity \cite{Sof}.

A natural question then follows: is it true that in general, solutions of dispersive equations converge in appropriate norm ($L^2$ or $H^1$) to a free wave and independently moving localized parts? In fact this is precisely the statement of asymptotic completeness in the case of $N$-body scattering \cite{D-Ger,Gr,H-Sig1,H-Sig2,SS1987,SSInvention}. In this case the possible outgoing clusters are clearly identified, as bound states of subsystems. But when the interaction term includes time dependent potentials (even localized in space) and more general nonlinear terms, we do not have an \emph{a priori} knowledge of the possible asymptotic states.

In fact, there are no general scattering results for localized time-dependent potentials. The exceptions are charge transfer Hamiltonians \cite{Yaj1, Gr1,Wul,Per, RSS,DSY}, decaying-in-time potentials and small potentials \cite{ How, RS}, time periodic potentials \cite{Yaj2, How} and random (in time) potentials \cite{BeSof}. See also \cite{BeSof1}. For potentials with asymptotic energy distribution more could be done \cite{SS1988}. A recent progress for more general localized potentials without smallness assumptions is obtained in \cite{SW2020}. Some tools from this work will be used in this paper.

Turning to the nonlinear case, Tao \cite{T1,T2,T3} has shown that the asymptotic decomposition holds for the nonlinear Schr\"odinger equation (NLS) with inter-critical nonlinearities in 3 or higher dimensions, in the case of radial initial data. In particular, in very high dimensions, and with an interaction that is a sum of a smooth compactly supported potential and a repulsive nonlinearity, Tao was able to show that the localized part is smooth and localized. In other cases, Tao showed the non-free part is only \emph{weakly localized} and smooth. Tao's work uses direct estimates of the incoming and outgoing parts of the solution to control the nonlinear part, via Duhamel representation. In a certain sense, it is in the spirit of Enss's work. See also \cite{Rod-Tao}.

A new approach due to Liu-Soffer \cite{Liu-Sof1,Liu-Sof2} is based on proving \emph{a priori} estimates on the full dynamics, which hold in suitably localized regions of the extended phase-space. In this way it was possible to show the asymptotic decomposition for general localized interactions, including time and space dependent interactions. Radial initial data is assumed in the nonlinear case. More detailed information is obtained on the localized part of the solution. Besides being smooth, its expanding part (if it exists) can grow at most like $|x|\leq \sqrt t,$ and furthermore, is concentrated in a thin set of the extended phase-space. The free part of the solution concentrates on the \emph{propagation set} where $x=vt, v=2p,$ and $p$ being the dual to the space variable, the momentum, is given by the operator $-i\nabla_x.$ The weakly localized part is found (in some sense) to be localized in the regions (see~\cite{Liu-Sof2}) where $$|x|/t^{\alpha}\sim 1  \quad\textbf{and} \quad |p|\sim  t^{-\alpha}, \quad\forall ~0<\alpha \leq 1/2.$$ It therefore predicts that the spreading part follows a self-similar pattern. The radial restriction was subsequently removed for Schr\"odinger equations by Soffer--Wu \cite{SW20221}, and an analogous framework for Klein--Gordon type equations was developed in \cite{SW20222}.

The question is therefore, do there exist solutions which are weakly localized but not localized? For equations which are not dispersive/hyperbolic many such solutions are known. But for hyperbolic/dispersive equations it is harder to see how such global solutions can emerge. In special cases such solutions were constructed \cite{MP2020,MZZ2003,Sch2023}. Such solutions have infinite $L^2$ norm and the equation is self-similar. In the energy critical nonlinear wave equation, the conformal symmetry implies a specific structure of scaling. Indeed in these models, it was possible to show that the asymptotic states also include self-similar solitons \cite{Duy-KMM,DKM2}. See also cited references in \cite{Duy-KMM,DKM2}. These solutions are global and preserve the energy (but have divergent $L^2$ norm). In \cite{S2021} it is shown that for small data there is a spreading weakly localized part for the KdV type equation.

In the case of NLS we have both energy and $L^2$ conservation, and no conformal symmetry. The aim of this work is to construct self-similar solutions for NLS type equations which are global and with non-zero mass ($L^2$ norm). We do this first for linear problems with time-dependent potentials that decay like $r^{-2}$ in space at infinity, and then use it to show that the addition of nonlinear terms does not change much the situation. We do not know yet if one can get such solutions for a purely nonlinear equation, but it is more plausible now. Our approach to this problem can be described in the language of scattering theory: We construct a channel of scattering where the asymptotic dynamics are given by a scaling transformation, which is unitary in $L^2,$ and a phase transformation. This is to be expected in view of the recent work \cite{Liu-Sof1,Liu-Sof2} showing that the asymptotic solutions (in the radial case) concentrate on thin sets of the phase-space, corresponding to self-similar solutions. This is the case for both the free wave and the weakly localized part.

\medskip

\subsection{Notation}\label{sec:notation}
We write $\s^p_x:=L^p(\mathbb{R}^n)$ and use the mixed notation $\s^\infty_t\s^q_x$ for $L^\infty_tL^q_x(\mathbb{R}^n\times [1,\infty))$. We write $X\lesssim Y$ to mean $X\leq C\,Y$ for some absolute constant $C>0$, and $X\lesssim_\alpha Y$ to mean $X\leq C(\alpha)Y$ with the constant allowed to depend on the parameter $\alpha$. For a real-valued function $A(t)$ and a positive-valued function $B(t)$, we write $A(t)\sim B(t)$ as $t\to\infty$ to mean that there exist constants $C_1,C_2>0$, $\tau\geq 1$, and a sign $\sigma\in\{+,-\}$ such that $\sigma C_1 B(t)\leq A(t)\leq \sigma C_2 B(t)$ for all $t\geq \tau$ (in particular, $A$ has constant sign on $[\tau,\infty)$). Following \cite{SW20221}, the symbols $F_c$ and $F_2$ (with an inequality as argument, e.g.\ $F_c(|x-2tp|/t^\alpha\leq 1)$, $F_2(|x|/\langle s\rangle^\beta\leq 1)$) always denote smooth, non-negative bump functions on $\mathbb{R}^n$ that equal $1$ on the region where the indicated inequality holds and vanish on its complement outside a thin transition layer; we write $\bar F_c:=1-F_c$ and $\bar F_2:=1-F_2$.

\subsection{Main results}
The aim of this subsection is to state the three main theorems of the paper. After describing the general equation \eqref{SE} and our standing assumptions on the self-similar potential and the scaling function $g(t)$, we proceed in three layers:
\begin{itemize}
%\item a self-similar dispersion estimate for the weakly localized part (Theorem \ref{wlem});
\item the existence of a non-trivial self-similar channel for the linear mass-critical system \eqref{linear} when $n\geq3$ (Theorem \ref{thm1});
\item a two-bubble extension to a linear mixed model where both non-trivial self-similar and localized asymptotics exist when $n\geq5$ (Theorem \ref{thm10});
\item a two-bubble result for a focusing nonlinear Schr\"odinger equation \eqref{NP}, where both non-trivial self-similar and localized asymptotics exist when $n\geq5$ (Theorem \ref{thm2}).
\end{itemize}

We consider the general class of nonlinear Schr\"odinger-type equations of the form
\begin{equation}
\begin{cases}
i\partial_t\psi(t) - H_0\psi(t) = V(x,t)\psi(t) + \mathcal{N}(|\psi(t)|)\psi(t) \\
\psi(x,1) = \psi_0 \in H^s(\mathbb{R}^n)
\end{cases}, \quad (x,t) \in \mathbb{R}^n \times [1,\infty) \label{SE}
\end{equation}
in space dimension $n\geq 3$, where $H_0 := -\Delta_x$ and $H^s(\mathbb{R}^n)$ is the $L^2$-Sobolev space in $x$. We take $s=0$ in the linear case and $s=1$ in the nonlinear case. Throughout the paper we only consider $t\geq 1$; in particular, the scaling function $g(t)$ is defined on $[1,\infty)$ and belongs to $C^2([1,\infty))$, and all time integrals, suprema, infima, and mixed-norm domains in the time variable are understood to be taken over $[1,\infty)$.

\subsubsection{Assumptions} We will make the following assumptions.
\begin{assumption}
The linear potential $V(x,t)$ takes the self-similar form: 
\begin{equation}
V(x,t) = g(t)^{-2} V\left( \frac{x}{g(t)} \right)
\end{equation}
where the time-independent potential $V(x)$ and the scaling factor $g(t)$ satisfy the following conditions:
\begin{enumerate}
    \item The stationary Hamiltonian $H := H_0 + V(x)$ has exactly one negative eigenvalue associated with a unique (up to a phase factor) normalized eigenfunction $\psi_b$.
    \item The scaling function $g\in C^2([1,\infty))$, and there exists $\epsilon \in (0, 1/2)$ such that the following symbol-type estimates hold for $k=0,1,2$:
    \begin{equation}
        \sup_{t \geq 1}
        \jBra{t}^{\epsilon + k}
        |g^{(k)}(t)| < \infty.
    \end{equation}
\end{enumerate}
Sharper, quantitative versions of these conditions, which will be used throughout the paper, are collected in \eqref{g(t)} and \eqref{Linear:con} below.
\end{assumption}

\subsubsection{The linear mass-critical system}
We first study the linear problem obtained from \eqref{SE} by setting $\mathcal{N}\equiv 0$ and choosing the initial datum to be a dilated bound state:
\eq
\begin{cases}
i\partial_t\psi(t)=H_0\psi(t)+g(t)^{-2} V(\frac{x}{g(t)})\psi(t) \\
\psi(x,t_0)=e^{- iD \ln g(t_0)}\psi_b(x)\in \s^2_x(\mathbb{R}^n)
\end{cases}\label{linear}, \quad n\geq 3
\eeq
for some $t_0\geq 1$ to be chosen later. The scaling function $g\in C^2([1,\infty))$ is required to satisfy the following sharper version of the symbol-type bounds of the standing Assumption: there exist constants $c_g\in (0,1)$ and $\epsilon\in (0,1/2)$ such that
\eq
\begin{cases}
\inf\limits_{t\geq 1}g(t)>0\\
 g(t)\sim \langle t\rangle^{\epsilon}\text{ as }t\to \infty \\
g(t)\sim tg'(t)\sim t^2g''(t)\text{ as }t\to \infty \\
\lim\limits_{t\to \infty}\frac{g(t)-2tg'(t)}{g(t)}= c_g
\end{cases}\label{g(t)}.
\eeq
Here and throughout the paper, the symbol $\sim$ stands for the two-sided bound defined in Section \ref{sec:notation} and we write $p:=-i\nabla_x$ for the momentum operator and $D:=\tfrac{1}{2}(x\cdot p+p\cdot x)$ for the generator of dilations on $\s^2_x$. The prototypical scaling function satisfying \eqref{g(t)} is $g(t)=\langle t\rangle^{\epsilon}$ with $\epsilon\in(0,1/2)$, for which a direct computation gives $c_g=1-2\epsilon\in(0,1)$; it is helpful for the reader to keep this concrete example in mind when parsing the general conditions above.

The time-independent potential $V(x)$ and the Hamiltonian $H:=H_0+V(x)$ are assumed to satisfy that 
\begin{equation}\label{eq: H lambda}
\parbox{0.9\linewidth}{$H$ has a unique normalized eigenstate (up to a phase factor) $\psi_b(x)$ with negative eigenvalue $\lambda<0$,}
\end{equation}
together with the decay conditions
\eq
\begin{cases}
\langle x\rangle V(x) \in \s^\infty_x\\
V(x)\in \s^2_x\cap \s^{n/2}_{\mathrm{loc}}(\mathbb{R}^n)
\end{cases}\label{Linear:con}
\eeq
on $V$. Throughout the paper we write $P_c$ for the spectral projection onto the continuous spectrum of $H$. Following \cite{RS}, a non-trivial function $\psi$ satisfying $H\psi=0$ in the sense of distributions and belonging to $\langle x\rangle^{\sigma}\s^2_x(\mathbb{R}^n)\setminus \s^2_x(\mathbb{R}^n)$ for some $\sigma>1/2$ is called a \emph{zero-energy resonance} of $H$; we say that \emph{$0$ is regular for $H$} when $0$ is neither an eigenvalue nor a resonance of $H$.

\begin{lemma}\label{lem:Hspectral}
Suppose $V$ satisfies \eqref{Linear:con} and $H=H_0+V$ satisfies \eqref{eq: H lambda}. Then
\begin{enumerate}
\item\label{it:0reg} $0$ is a regular point of $H$ (i.e.\ neither an eigenvalue nor a resonance);
\item\label{it:Dpsib} $\langle x\rangle D\psi_b\in \s^2_x$ and $\|D\psi_b\|_{\s^\infty_x}\lesssim 1$;
\item\label{it:res} $\bigl\|(\lambda-H)^{-1}P_c\bigr\|_{\s^\infty_x\to \s^\infty_x}\lesssim_\lambda 1$.
\end{enumerate}
\end{lemma}
\begin{proof}
\emph{(\ref{it:0reg})}. The unique-negative-eigenvalue hypothesis rules out $0$ being an eigenvalue of $H$. Under \eqref{Linear:con}, the absence of positive embedded eigenvalues (and of a zero resonance) follows from Ionescu--Jerison \cite{IJ2003}; for $n=3$ classical results are also available in Komech--Kopylova \cite{KK2012}, Section~4.

\emph{(\ref{it:Dpsib})}. From $H\psi_b=\lambda\psi_b$ with $\lambda<0$ one obtains
\[
\psi_b=(H_0-\lambda)^{-1}V\psi_b.
\]
Since $\langle x\rangle V\in \s^\infty_x$ and $(H_0-\lambda)^{-1}$ has an exponentially-decaying kernel at the negative energy $\lambda$, a standard bootstrap yields exponential decay of $\psi_b$ and of its derivatives. Consequently, $\langle x\rangle D\psi_b\in \s^2_x$, and since $D\psi_b=-ix\cdot\nabla\psi_b-\tfrac{in}{2}\psi_b$, the exponential decay of $\psi_b$ and $\nabla\psi_b$ also yields $\|D\psi_b\|_{\s^\infty_x}\lesssim 1$.

\emph{(\ref{it:res})}. Under \eqref{Linear:con}, the wave operators $W_\pm:=\mathrm{s}\text{-}\lim_{t\to\pm\infty}e^{itH}e^{-itH_0}$ exist and are bounded on $\s^p_x$ for all $p\in[1,\infty]$ after a high-frequency cut-off (Soffer--Wu \cite{SW2020}); together with asymptotic completeness (which holds in our setting since $0$ is regular for $H$ by part \eqref{it:0reg}), they intertwine the absolutely continuous parts of $H_0$ and $H$ as
\[
(\lambda-H)^{-1}P_c=W_+(\lambda-H_0)^{-1}W_+^*.
\]
Since $\lambda<0$, the operator $(\lambda-H_0)^{-1}=(-\Delta-\lambda)^{-1}$ is bounded on $\s^\infty_x$; moreover, its symbol $1/(|p|^2-\lambda)$ has rapid decay at high frequency, so the high-frequency cut-off needed for the $\s^\infty_x$-boundedness of $W_\pm,W_\pm^*$ in \cite{SW2020} is automatic in this composition. Combining the three bounds yields $\|(\lambda-H)^{-1}P_c\|_{\s^\infty_x\to \s^\infty_x}\lesssim_\lambda 1$.
\end{proof}

We refer to system \eqref{linear} as the \emph{mass-critical system}: the time-dependent potential $g(t)^{-2}V(x/g(t))$ has the same $\s^2_x$ scaling as $H_0$, so the linear and (later) nonlinear terms compete on the same scale.

\subsubsection{Self-similar decomposition and a rescaled time variable}
We split $\psi(t)$ along the one-dimensional self-similar direction spanned by the dilated bound state:
\eq
\psi(t)=\tilde a(t)e^{- iD \ln (g(t)) }\psi_b(x)\oplus \psi_c(t)\label{decomp}
\eeq
with the orthogonality condition
\eq
( e^{- iD \ln (g(t)) }\psi_b(x),\psi_c(t))_{\s^2_x}=0,\label{orth}
\eeq
so that, since $e^{-iD\ln g(t)}$ is unitary on $\s^2_x$,
\eq
\tilde a(t):=(e^{-iD\ln g(t)}\psi_b(x),\psi(t))_{\s^2_x}=(\psi_b(x),e^{iD\ln g(t)}\psi(t))_{\s^2_x}\label{a}
\eeq
is the $\s^2_x$-overlap of $\psi(t)$ with the moving self-similar profile $e^{-iD\ln g(t)}\psi_b(x)$, and $\psi_c(t)$ is the orthogonal remainder. Showing that this overlap does not vanish as $t\to\infty$ is the content of our main linear result.

To track $\tilde a(t)$, it is convenient to pass to a rescaled time. Define $T: [1,\infty)\to [0,\infty)$ by $t\mapsto s=\int_1^t g(u)^{-2}\,du$; since $T(1)=0$ and $T$ is strictly increasing, the inverse $T^{-1}: [0,\infty)\to [1,\infty)$ is well defined. By \eqref{g(t)},
\eq
s=T(t)=\int_1^t g(u)^{-2}\,du\sim \int_1^t \langle u\rangle^{-2\epsilon}\,du\sim \langle t\rangle^{1-2\epsilon}\to \infty \text{ as }t\to \infty.
\eeq
After a unitary conjugation by $e^{iD\ln g(t)}$ and passage to the rescaled time $s$, the time-dependent operator $H_0+g(t)^{-2}V(x/g(t))$ is transformed into the autonomous Hamiltonian $H$ plus a slowly-varying (but non-integrable) dilation perturbation $f(s)D$, with $f(s)\sim\langle s\rangle^{-1}$ as in~\eqref{fs}. This recasts the persistence of $\tilde a(t)$ as an ionization problem for $H$; the reduction is carried out in Section~\ref{outline}.

\subsubsection{Main results for the linear problem}
Our first main result asserts that $\psi(t)$ of \eqref{linear} contains a non-trivial self-similar component
\eq
\psi_{w,b}(t):=\tilde a(t)e^{-iD\ln(g(t) )}\psi_b(x)
\eeq
in the following sense:
\begin{enumerate}
\item \label{item:A}For every solution $\psi(t)$ of \eqref{linear}, with $\tilde a(t)$ as in \eqref{a}, the limit
\eq
\tilde{A}(\infty):=\lim\limits_{t\to \infty} e^{i\lambda T(t)}\tilde{a}(t)\label{tildea}
\eeq
exists.
\item \label{item:B}For all $t_0>1$ sufficiently large, the particular initial condition
\eq
\psi(t_0)=e^{-i D\ln(g(t_0) )}\psi_b(x)
\eeq
prescribed in \eqref{linear} yields
\eq
|\tilde{A}(\infty)| >0,
\eeq
which implies
\eq\liminf\limits_{t\to \infty}|\tilde{a}(t)|\gtrsim 1.\label{eq10}
\eeq
\end{enumerate}
The state $e^{-iD\ln g(t)}\psi_b(x)$ is the self-similar profile
\eq
e^{-i D\ln(g(t) )}\psi_b(x)=\frac{1}{g(t)^{n/2}}\psi_b(x/g(t)),
\eeq
which, since $e^{-iD\ln g(t)}$ is a unitary operator on $\s^2_x$, has the same $\s^2_x$-norm as $\psi_b$:
\eq
\|e^{-i D\ln(g(t) )}\psi_b(x)\|_{\s^2_x(\mathbb{R}^n)}=\|\psi_b(x)\|_{\s^2_x(\mathbb{R}^n)}=1.
\eeq
Thus item~(\ref{item:A}) above asserts that, after compensating for the bound-state oscillation via the factor $e^{i\lambda T(t)}$, the $\s^2_x$-overlap of $\psi(t)$ with the moving self-similar profile has a definite limit; item~(\ref{item:B}) says that, for the specific initial datum prescribed in \eqref{linear}, this limit is non-trivial. The analytic core of this statement is isolated as Lemma \ref{lem:wself4} in Section \ref{sec:toolbox}.

\begin{theorem}\label{thm1}Let $\tilde a(t)$ be as in \eqref{a}, and let $\psi(t)$ be the solution of \eqref{linear}. Suppose $V$ satisfies \eqref{Linear:con}, $H=H_0+V$ satisfies \eqref{eq: H lambda}, and $g$ satisfies \eqref{g(t)}. Then for $n\geq 3$,
\eq\label{eq: tA exist}
\tilde{A}(\infty):=\lim\limits_{t\to \infty} e^{i\lambda T(t)}\tilde a(t)
\eeq
exists, the decomposition \eqref{decomp}--\eqref{a} of $\psi(x,t)$ holds with $\tilde{a}(t)$ satisfying \eqref{eq10}, and the $g(t)$-self-similar channel wave operator
\eq
\Omega_{g}^*\psi(0):=w\text{-}\lim\limits_{s\to \infty}e^{isH}e^{ iD\ln (g( T^{-1}(s)))}\psi(T^{-1}(s))
\eeq
exists in $\s^2_x$ with
\eq
\Omega_{g}^*\psi(0)=\tilde{A}(\infty)\psi_b(x).\label{Omegag}
\eeq
\end{theorem}
\subsubsection{The mixed model: an additional fixed bound state}
Building on \eqref{linear}, we now add a second, time-independent potential $W(x)$ such that $H_0+W$ admits a bound state, and consider the following mixed model:
\eq
\begin{cases}
i\partial_t\psi(t)=H_0\psi(t)+W(x)\psi(t)+g(t)^{-2}V(\frac{x}{g(t)})\psi(t)\\
\psi(x,t_0)=\psi_d(x)+e^{-i D\ln(g(t_0) )}\psi_b(x)\in H^1(\mathbb{R}^n)
\end{cases}\quad (x,t)\in \mathbb{R}^n\times [1,\infty)\label{NP0},
\eeq
when $n\geq 5$, where $H_0+W$ is assumed to satisfy
\begin{equation}\label{eq: W lambda0}
\parbox{0.9\linewidth}{$H_0+W$ has a unique normalized eigenstate (up to a phase factor) $\psi_d(x)$ with negative eigenvalue $\lambda_0<0$,}
\end{equation}
together with the decay conditions
\eq
\begin{cases}
W(x)\in \s^2_x(\mathbb{R}^n)\cap \s^\infty_x(\mathbb{R}^n)\\
x\cdot\nabla V(x)\in \s^\infty_x(\mathbb{R}^n)
\end{cases}\label{W}
\eeq
on $(W,V)$.
Conditions~\eqref{Linear:con}, \eqref{eq: W lambda0} and~\eqref{W} yield the self-adjointness of $H_0+W+g(t)^{-2}V(\frac{x}{g(t)})$ on $H^2(\mathbb R^n)$ by Kato's theorem. The initial datum is taken to be the sum of the fixed bound state $\psi_d(x)$ of $H_0+W$ and the dilated bound state $e^{-iD\ln g(t_0)}\psi_b(x)$ of the self-similar trap. If $\psi(t_0)\in H^1(\mathbb{R}^n)$, the solution of \eqref{NP0} satisfies the following \emph{a priori} $H^1$ bound.
\begin{lemma}[\emph{A priori} $H^1$ bound for \eqref{NP0}]\label{lem:NP0H1}
Suppose $V$ satisfies \eqref{Linear:con}, the pair $(V,W)$ satisfies \eqref{eq: W lambda0} and \eqref{W}, and $g$ satisfies \eqref{g(t)}. Then every solution $\psi\in C\bigl([t_0,\infty);\,H^1(\mathbb{R}^n)\bigr)$ of \eqref{NP0} satisfies the \emph{a priori} bound
\eq
\sup_{t\geq t_0}\|\psi(t)\|_{H^1}\leq C_0,\label{NP0:H1bound}
\eeq
in which the constant $C_0$ depends only on $\|\psi(t_0)\|_{H^1}$, $\|W\|_{\s^\infty_x}$, $\|V\|_{\s^\infty_x}$, $\|x\cdot\nabla V\|_{\s^\infty_x}$, and $g(t_0)$.
\end{lemma}
\begin{proof}
Since $H_0+W+g(t)^{-2}V(\frac{x}{g(t)})$ is self-adjoint, the $\s^2_x$-mass is conserved:
\eq
\|\psi(t)\|_{\s^2_x}=\|\psi(t_0)\|_{\s^2_x},\qquad t\geq t_0.
\eeq
Set $\tilde V(x,t):=g(t)^{-2}V(x/g(t))$ and define the energy
\eq
E(t):=\bigl(\psi(t),(H_0+W+\tilde V(t))\psi(t)\bigr)_{\s^2_x}.
\eeq
Since $H(t):=H_0+W+\tilde V(t)$ is self-adjoint and $H_0$, $W$ are time-independent,
\eq
\tfrac{d}{dt}E(t)=\bigl(\psi(t),(\partial_t\tilde V)(t)\psi(t)\bigr)_{\s^2_x}.
\eeq
A direct computation gives
\eq
\partial_t\tilde V(x,t)=-\frac{g'(t)}{g(t)^3}\bigl[2V(y)+y\cdot\nabla V(y)\bigr]_{y=x/g(t)},
\eeq
so that $\|\partial_t\tilde V(\cdot,t)\|_{\s^\infty_x}\leq \frac{|g'(t)|}{g(t)^3}\bigl(2\|V\|_{\s^\infty_x}+\|x\cdot\nabla V\|_{\s^\infty_x}\bigr)$. Combining with $\s^2_x$-conservation,
\eq
\Bigl|\tfrac{d}{dt}E(t)\Bigr|\leq \frac{|g'(t)|}{g(t)^3}\bigl(2\|V\|_{\s^\infty_x}+\|x\cdot\nabla V\|_{\s^\infty_x}\bigr)\|\psi(t_0)\|_{\s^2_x}^2.\label{dEdt}
\eeq
By \eqref{g(t)} one has $g(t)\sim\langle t\rangle^\epsilon$ with $\epsilon\in(0,1/2)$, so $\int_{t_0}^\infty|g'|/g^3\,ds<\infty$; integrating \eqref{dEdt} gives $\sup_{t\geq t_0}|E(t)|<\infty$. Since $H_0=-\Delta_x$,
\eq
\|\nabla\psi(t)\|_{\s^2_x}^2=E(t)-(\psi(t),W\psi(t))-(\psi(t),\tilde V(t)\psi(t))\leq E(t)+\bigl(\|W\|_{\s^\infty_x}+g(t)^{-2}\|V\|_{\s^\infty_x}\bigr)\|\psi(t_0)\|_{\s^2_x}^2,
\eeq
which, together with $\s^2_x$-conservation, yields \eqref{NP0:H1bound}. The initial energy $E(t_0)$ is controlled by $\|\psi(t_0)\|_{H^1}^2+(\|W\|_{\s^\infty_x}+g(t_0)^{-2}\|V\|_{\s^\infty_x})\|\psi(t_0)\|_{\s^2_x}^2$.
\end{proof}

Before stating the theorem, we fix the terminology that will be used throughout the paper. By a \emph{bubble} in a solution $\psi(t)$ we mean a non-trivial, asymptotically persistent coherent profile --- that is, either a spatially fixed localized shape (a \emph{localized bubble}, e.g.\ a bound state of $H_0+W$ or a soliton) or a profile moved and/or dilated along an explicit family of symmetries (a \emph{self-similar bubble}, e.g.\ the dilated bound state $e^{-iD\ln g(t)}\psi_b(x)$ of Theorem \ref{thm1}). The next theorem shows that the solution of the mixed model above asymptotically carries at least \emph{two bubbles}: a non-trivial self-similar bubble (as in Theorem \ref{thm1}) plus a non-trivial localized bubble near the origin, persisting in the direction of $\psi_d$.
\begin{theorem}\label{thm10}Let $\tilde a(t)$ be as in \eqref{a}, and let $\psi(t)$ be the solution of \eqref{NP0}. Suppose $V$ satisfies \eqref{Linear:con}, $H=H_0+V$ satisfies \eqref{eq: H lambda}, $H_0+W$ satisfies \eqref{eq: W lambda0}, the pair $(V,W)$ satisfies \eqref{W}, and $g$ satisfies \eqref{g(t)}. Then, for $n\geq 5$, $\epsilon\in(2/n, 1/2)$, and $t_0\geq 1$ sufficiently large,
\eq
\tilde{A}(\infty):=\lim\limits_{t\to \infty} e^{i\lambda T(t)}\tilde a(t)\label{tildeAthm10}
\eeq
exists, with $\tilde a(t)$ satisfying \eqref{eq10}, and there exists a constant $c_d>0$ such that
\eq
\liminf_{t\to \infty} | ( \psi(x,t),\psi_d(x))_{\s^2_x}|\geq c_d.\label{cd0}
\eeq
Moreover, the $g(t)$-self-similar channel wave operator
\eq
\Omega_{g}^*\psi(0):=w\text{-}\lim\limits_{s\to \infty}e^{isH}e^{ iD\ln (g( T^{-1}(s)))}\psi(T^{-1}(s))
\eeq
exists in $\s^2_x$ and
\eq
\Omega_{g}^*\psi(0)=\tilde{A}(\infty)\psi_b(x).\label{Omegag10}
\eeq
\end{theorem}

The two bubbles delivered by Theorem \ref{thm10} are of different natures. The self-similar bubble $\tilde a(t)e^{-iD\ln g(t)}\psi_b(x)$ is a \emph{directional} statement: a definite portion of the $\s^2_x$ mass stays aligned with the moving profile $e^{-iD\ln g(t)}\psi_b(x)$. The localized bubble, on the other hand, is captured by \eqref{cd0}, which says that the solution $\psi(x,t)$ retains a non-trivial projection on the fixed bound state $\psi_d$; this is strictly stronger than merely keeping $\s^2_x$ mass inside a fixed ball. Theorem \ref{thm2} below will replace this directional statement by the weaker, but still non-trivial, ball-mass estimate \eqref{con100}.

\subsubsection{Application: a focusing nonlinear Schr\"odinger equation}
The mixed-model analysis is robust enough to absorb a focusing nonlinearity: the role of the fixed bound state $\psi_d$ is played by a soliton of an associated stationary nonlinear problem. Paralleling the mixed model of Theorem \ref{thm10}, we therefore consider the focusing nonlinear Schr\"odinger equation
\eq
\begin{cases}
i\partial_t\psi(t)=H_0\psi(t)+g(t)^{-2}V(\frac{x}{g(t)})\psi(t)+\mathcal{N}(|\psi(t)|)\psi(t)\\
\psi(x,t_0)=\psi_s(x)+e^{-i D\ln(g( t_0) )}\psi_b(x)\in H^1(\mathbb{R}^n)
\end{cases},\quad (x,t)\in \mathbb{R}^n\times [1,\infty)\label{NP},
\eeq
when $n\geq 5$, in the focusing radial regime: we impose $\mathcal{N}\leq 0$ and require that both $\psi(t_0)$ and $V(x)$ be radial in $x$. In addition, assume $V$ satisfies \eqref{Linear:con} and $H:=H_0+V(x)$ satisfies \eqref{eq: H lambda} (so Lemma \ref{lem:Hspectral} applies, providing the bound state $\psi_b$), and assume $\psi_s(x)$ is a soliton of
\eq
i\partial_t\phi=H_0\phi+\mathcal{N}_{F,0}(|\phi|)\phi,\label{NF00}
\eeq
i.e., $\psi_s$ satisfies the stationary equation
\eq
(H_0+2\mathcal{N}_{F,0}(|\psi_s(x)|))\psi_s(x)=E\psi_s(x), \text{ for some }E<0,\label{soliton}
\eeq
where
\eq
\mathcal{N}_{F,0}(k):=\int_0^k dq q \mathcal{N}(q)/k^2<0.
\eeq
As $t_0\to\infty$ the (localized) soliton $\psi_s$ and the dilated bound state $e^{-iD\ln g(t_0)}\psi_b(x)=g(t_0)^{-n/2}\psi_b(x/g(t_0))$ separate in physical space (the latter spreading at scale $g(t_0)$), so the bilinear form of $\psi(t_0)=\psi_s+e^{-iD\ln g(t_0)}\psi_b$ converges to that of $\psi_s$ alone, while $g(t_0)^{-2}V(x/g(t_0))$ vanishes uniformly. Combined with \eqref{soliton} and $E<0$, this yields, for all $t_0\geq \tau$ sufficiently large,
\eq
(\psi(t_0),(H_0+\frac{1}{g(t_0)^{2}}V(\frac{x}{g(t_0)})+2\mathcal{N}_{F,0}(|\psi(t_0)|))\psi(t_0)  )_{\s^2_x}\leq \frac{E}{2}\|\psi_s(x)\|_{\s^2_x}^2,\label{Rt0}
\eeq
with $\psi(t_0)=\psi_s(x)+e^{-i D\ln(g(t_0))}\psi_b(x)$. We further assume the nonlinearity satisfies the polynomial bound
\eq
| \mathcal{N}(k)|\lesssim |k|^\beta,\text{ for some }\beta>0.\label{NFandN}
\eeq
Under \eqref{soliton} and \eqref{NFandN}, the solution of \eqref{NP} carries two bubbles --- one self-similar and one localized near the origin --- in the sense made precise by Theorem \ref{thm2} below. The $H^1$ bound
\eq
\sup_{t\geq t_0}\|\psi(t)\|_{H^1}\lesssim_{\psi(t_0)} 1\label{NP:H1bound}
\eeq
is a consequence of an asymptotic-energy argument, which we record as Lemma \ref{lem:NPH1} below; for the concrete focusing saturated nonlinearity treated in Section \ref{sectionNP}, \eqref{NP:H1bound} is verified by a direct computation in Lemma \ref{typicalNP}.

The radial assumption, combined with $\psi(t)\in H^1(\mathbb R^n)$, implies the classical Strauss radial decay (\cite{Strauss1977}, Radial Lemma 1)
\eq
|\psi(x,t)|\lesssim \frac{1}{|x|^{\frac{n-1}{2}}}, \quad |x|\geq 1.\label{rad}
\eeq
As a consequence of \eqref{NFandN} together with the Sobolev embedding $H^1\hookrightarrow L^p$, every $f\in H^1(\mathbb{R}^n)$ satisfies
\eq
\begin{cases}
\|\mathcal{N}(|f(x)|)\|_{\s^2_x}\leq C(\|f(x)\|_{H^1})\\
\| \mathcal{N}(|f(x)|)f(x)\|_{\s^2_x}\leq C(\|f(x)\|_{H^1})\label{con:N}
\end{cases},
\eeq
which we shall use freely below. Combining \eqref{NFandN} with \eqref{rad} yields, for all $M\geq 1$, $t\geq t_0$ and $t_0$ sufficiently large,
\eq
\left|(\psi(t), \mathcal{N}_{F,0}(|\psi(t)|)\psi(t))_{\s^2_x}\right|\leq  \frac{1}{M^{\frac{(n-1)\beta}{2}}}C(\| \psi(t)\|_{H^1})+\| \chi(|x|\leq M)\psi(t)\|_{\s^2_x}C(\| \psi(t)\|_{H^1}).\label{ACEN}
\eeq

\begin{lemma}[Global $H^1$ solution and \emph{a priori} bound for \eqref{NP}]\label{lem:NPH1}
Suppose $V$ satisfies \eqref{Linear:con}, $H=H_0+V$ satisfies \eqref{eq: H lambda}, $g$ satisfies \eqref{g(t)}, and $\mathcal N$ satisfies \eqref{soliton} and \eqref{NFandN}. Assume in addition that $x\cdot\nabla V\in \s^\infty_x(\mathbb{R}^n)$, $\mathcal{N}\in \s^\infty_x(\mathbb{R}_+)$, and $|\mathcal{N}_F(k)|\lesssim k^2$, where $\mathcal{N}_F(k):=\int_0^k q^2\mathcal{N}'(q)\,dq$. Then, for $t_0\geq 1$ sufficiently large, equation \eqref{NP} admits a global solution $\psi\in C([t_0,\infty);H^1(\mathbb{R}^n))$ satisfying \eqref{NP:H1bound}.
\end{lemma}
\begin{proof}
Local well-posedness in $H^1$ follows from a standard contraction argument: since \eqref{Linear:con} gives $V\in \s^\infty_x$ and $\|g(t)^{-2}V(\cdot/g(t))\|_{\s^\infty_x}\leq g(t_0)^{-2}\|V\|_{\s^\infty_x}$ uniformly in $t\geq t_0$, the linear evolution $U(t,s)$ generated by $H_0+g(t)^{-2}V(x/g(t))$ is unitary on $\s^2_x$ and satisfies $\|U(t,s)\|_{H^1\to H^1}\leq e^{C|t-s|}$ on bounded intervals; together with $\mathcal N\in \s^\infty_x$ (Lipschitz control of the nonlinear term in $H^1$), this yields a unique solution $\psi$ on a maximal interval $[t_0,t_*)$. The \emph{a priori} bound below shows $t_*=\infty$.

Set $\tilde V(x,t):=g(t)^{-2}V(x/g(t))$ and observe that, since $V,\mathcal{N}$ are real-valued, $\s^2_x$-mass is conserved: $\|\psi(t)\|_{\s^2_x}=\|\psi(t_0)\|_{\s^2_x}$ for all $t\geq t_0$. Using the chain rule, the total instantaneous energy
\eq
\mathcal E(t):=\bigl(\psi(t),\bigl(H_0+\tilde V(t)+\mathcal{N}(|\psi(t)|)\bigr)\psi(t)\bigr)_{\s^2_x}
\eeq
satisfies
\eq
\frac{d}{dt}\mathcal E(t)=(\psi(t),(\partial_t\tilde V)\psi(t))_{\s^2_x}+\frac{d}{dt}\!\int \mathcal{N}_F(|\psi(t)|)\,dx,
\eeq
since $\langle\psi,\partial_t[\mathcal{N}(|\psi|)]\psi\rangle=\partial_t\!\int\mathcal{N}_F(|\psi|)\,dx$ by the identity $k^2\mathcal{N}'(k)=\mathcal{N}_F'(k)$. The direct computation $\partial_t\tilde V(x,t)=-\frac{g'(t)}{g(t)^3}[2V(y)+y\cdot\nabla V(y)]|_{y=x/g(t)}$ and the assumption $x\cdot\nabla V\in \s^\infty_x$ give
\eq
\|\partial_t\tilde V(\cdot,t)\|_{\s^\infty_x}\leq \frac{|g'(t)|}{g(t)^3}\bigl(2\|V\|_{\s^\infty_x}+\|x\cdot\nabla V\|_{\s^\infty_x}\bigr),
\eeq
and \eqref{g(t)} yields $\int_{t_0}^\infty|g'|/g^3\,ds<\infty$, so $\int_{t_0}^\tau |\langle\psi,\partial_t\tilde V\,\psi\rangle|\,ds\lesssim \|\psi(t_0)\|_{\s^2_x}^2$ uniformly in $\tau$. Combining with $|\mathcal{N}_F(k)|\lesssim k^2$ and $\s^2_x$-conservation,
\eq
\Bigl|\int_{t_0}^\tau\!\frac{d}{dt}\!\int\mathcal{N}_F(|\psi|)\,dx\,dt\Bigr|=\Bigl|\int\mathcal{N}_F(|\psi(\tau)|)\,dx-\int\mathcal{N}_F(|\psi(t_0)|)\,dx\Bigr|\lesssim \|\psi(t_0)\|_{\s^2_x}^2,
\eeq
so $|\mathcal E(\tau)-\mathcal E(t_0)|\lesssim \|\psi(t_0)\|_{\s^2_x}^2$. Expanding
$$\mathcal E(\tau)=\|\nabla\psi(\tau)\|_{\s^2_x}^2+(\psi(\tau),\tilde V(\tau)\psi(\tau))+(\psi(\tau),\mathcal{N}(|\psi(\tau)|)\psi(\tau))$$
and using $\|\tilde V(\tau)\|_{\s^\infty_x}\leq g(\tau)^{-2}\|V\|_{\s^\infty_x}$ together with $\mathcal{N}\in \s^\infty_x$,
\eq
\|\nabla\psi(\tau)\|_{\s^2_x}^2\leq |\mathcal E(t_0)|+C\|\psi(t_0)\|_{\s^2_x}^2\bigl(1+\|V\|_{\s^\infty_x}+\|\mathcal{N}\|_{\s^\infty_x}\bigr),
\eeq
and the initial energy $\mathcal E(t_0)$ is controlled by $\|\psi(t_0)\|_{H^1}^2+C_0\|\psi(t_0)\|_{\s^2_x}^2$. This yields \eqref{NP:H1bound}.
\end{proof}

\begin{theorem}\label{thm2}Let $\tilde a(t)$ be as in \eqref{a}, and let $\psi(t)$ be the solution of \eqref{NP}. Suppose $V$ satisfies \eqref{Linear:con} and the $V$-part of \eqref{W} (i.e., $x\cdot\nabla V\in \s^\infty_x$), and $H=H_0+V$ satisfies \eqref{eq: H lambda} (so that Lemma \ref{lem:Hspectral} applies); that $V$ and $\psi(t_0)$ are radial in $x$; that $\mathcal{N}\leq 0$ satisfies \eqref{soliton}, \eqref{NFandN}, $\mathcal N\in \s^\infty_x(\mathbb R_+)$, and $|\mathcal N_F(k)|\lesssim k^2$ with $\mathcal N_F(k):=\int_0^k q^2\mathcal N'(q)\,dq$ (so that Lemma~\ref{lem:NPH1} applies); and that $g$ satisfies \eqref{g(t)}. Then, for $n\geq 5$, $\epsilon \in (2/n,1/2)$, and $t_0\geq 1$ sufficiently large,
\eq
\tilde{A}(\infty):=\lim\limits_{t\to \infty} e^{i\lambda T(t)}\tilde a(t)
\eeq
exists, the decomposition \eqref{decomp}--\eqref{a} of $\psi(x,t)$ holds with $\tilde a(t)$ satisfying \eqref{eq10}, and there exists some large number $M\geq 1$ such that
\eq
\liminf\limits_{t\to \infty} \|\chi(|x|\leq M) \psi(t)\|_{\s^2_x}\geq c'\label{con100}
\eeq
for some $c'>0$. Moreover, based on \eqref{wself4}, the $g(t)$-self-similar channel wave operator
\eq
\Omega_{g}^*\psi(0):=w\text{-}\lim\limits_{s\to \infty}e^{isH}e^{ iD\ln (g( T^{-1}(s)))}\psi(T^{-1}(s))
\eeq
exists in $\s^2_x$ and
\eq
\Omega_{g}^*\psi(0)=\tilde{A}(\infty)\psi_b(x).
\eeq
\end{theorem}
A concrete instance to which Theorem \ref{thm2} applies is the focusing saturated nonlinearity $\mathcal{N}(|\psi|)=-\lambda|\psi|/(1+|\psi|^2)$ with $g(t)=\langle t\rangle^\epsilon$, $\epsilon\in(2/5,1/2)$ and $\lambda>0$ sufficiently large, in $n=5$ space dimensions; see Lemma \ref{typicalNP} in Section \ref{sectionNP}.
\subsection{Outline of the proof}\label{outline}
We outline the proof in three blocks: the reduction used for the linear problem (Theorem \ref{thm1}), the adjustments needed to accommodate $W(x)$ or the nonlinearity $\mathcal N$ (Theorems \ref{thm10}, \ref{thm2}), and the energy argument that produces the second, localized bubble.

\subsubsection{Linear case: reduction to an ionization problem}

\paragraph{Step 1. Unitary conjugation.} The operator $D$ generates $\s^2_x$-preserving dilations, and a direct computation gives $e^{iD\ln g(t)}H_0e^{-iD\ln g(t)}=g(t)^{-2} H_0$. Conjugation by $e^{iD\ln g(t)}$ thus rescales $H_0$ and pulls the self-similar potential $g(t)^{-2}V(x/g(t))$ back to the time-independent $V(x)$. To be precise, setting
\eq
\tphi(t):=e^{i D\ln(g(t) )}\psi(x,t),\label{tphi}
\eeq
a direct computation (Lemma \ref{lem:tphi}) gives
\eq
\begin{cases}
i\partial_t\tphi=g(t)^{-2}H\tphi-(\partial_t[g(t)]g(t)^{-1})D\tphi\\
\tphi(t_0)=\psi_b(x)
\end{cases}\label{Linear:1},
\eeq
so that all the time dependence of $H_0+g(t)^{-2}V(x/g(t))$ is now concentrated in the scalar factor $g(t)^{-2}$ and in the adiabatic drift $-(\partial_t g/g)D$.

\paragraph{Step 2. Time change.} To absorb the $g(t)^{-2}$ prefactor we pass to the rescaled time $s=T(t)=\int_1^t g(u)^{-2}\,du$ and set
\eq
\phi(s):=\tphi(T^{-1}(s)),\quad t_0:=T^{-1}(s_0).
\eeq
Then \eqref{Linear:1} becomes
\eq
\begin{cases}
i\partial_s\phi=H\phi+f(s)D\phi\\
\phi(s_0)=\psi_b(x)
\end{cases}\label{Linear:2},
\eeq
with
\eq
f(s):=-(\partial_t[g(t)]g(t))\vert_{t=T^{-1}(s)}.\label{fs}
\eeq
By \eqref{g(t)} (see Lemma \ref{lem:fasymp}),
\eq
f(s)\sim \frac{1}{\langle s\rangle},\qquad |f'(s)|\lesssim \frac{1}{\langle s\rangle^2},\label{fsa}
\eeq
so $f(s)D$ is a slowly-varying ($f'(s)$ integrable) but \emph{non-integrable} time-dependent perturbation of the autonomous evolution generated by the time-independent Hamiltonian $H$; this borderline regime is what makes the oscillatory integration-by-parts in Step~4 below necessary.

\paragraph{Step 3. Ionization reduction.} Equation \eqref{Linear:2} is similar to the form studied in the ionization literature \cite{SW1}: an autonomous Hamiltonian $H$ with bound state $\psi_b$ perturbed by a slowly-decaying time-dependent term. The asymptotic behavior of the self-similar projection $\tilde a(t)$ from \eqref{a} is governed by the bound-state amplitude
\eq
a(s):=(\psi_b, \phi(s))_{\s^2_x},\label{defa}
\eeq
expressed in the rescaled time variable; the two are the same quantity with $\tilde a(t)=a(T(t))$. The linear result thus reduces to understanding the behavior of $a(s)$ as $s\to\infty$. The behavior obtained here should be contrasted with the Fermi-Golden-Rule ionization result of Soffer--Weinstein \cite{SW1}: for a generic time-dependent perturbation of the form $W(t,x)$ with suitable resonance between the bound-state energy $\lambda$ and the continuous spectrum, $|a(s)|\to 0$. By contrast, the self-similar perturbation $f(s)D$ considered here produces a non-resonant coupling between $\psi_b$ and the continuum, and $a(s)$ is shown to converge to a definite \emph{non-zero} limit. The novelty of Theorems \ref{thm1}--\ref{thm2} lies in identifying this persistence phenomenon and using it to construct a non-trivial self-similar asymptotic channel.

Both the convergence of $a(s)$ and its non-trivial lower bound make essential use of the spectral gap $\lambda<0= \inf\sigma_{\mathrm{ac}}(H)$ and of the assumption that $0$ is a regular point of $H$ (no threshold resonance or eigenvalue). Together these guarantee the local decay of $P_ce^{-isH}$, which is the engine powering the oscillatory integration-by-parts used in Step~4 below.

\paragraph{Step 4. Existence of the limit and a non-trivial lower bound.} The bound-state oscillation $e^{-i\lambda s}$ is removed by setting
\eq
A(s):=e^{i\lambda s }a(s),
\eeq
so that, from \eqref{Linear:2},
\eq
i\partial_s A(s)=e^{i\lambda s}f(s)\,(\psi_b, D P_c\phi(s))_{\s^2_x}.\label{dA}
\eeq
The key point is that although $f(s)\sim 1/\langle s\rangle$ is \emph{not} integrable, the right-hand side of \eqref{dA} is an oscillatory integral: integrating by parts in $s$, moving the factor $e^{i\lambda s}/(i\lambda)$ onto $f(s)(\psi_b, DP_c\phi(s))$ and using $f'(s)\sim 1/\langle s\rangle^2$ together with the continuous-spectrum dispersion of $e^{-isH}$, we obtain an $L^1_s$-bounded remainder. This yields the convergence of $A(s)$ and hence
\eq
A(\infty):=\lim\limits_{s\to \infty}e^{i\lambda s }a(s).
\eeq
From $\tilde{a}(T^{-1}(s))=a(s)$ we deduce
\eq
\tilde{A}(\infty)=A(\infty),
\eeq
which is~\eqref{eq: tA exist}. For the specific initial datum $\phi(s_0)=\psi_b$ in \eqref{Linear:2} one has $a(s_0)=1$; taking $s_0$ large enough so that the integration-by-parts tail is $<1/2$,
\eq
|\tilde A(\infty)|\geq \frac{1}{2}>0,
\eeq
which yields~\eqref{eq10}.

\subsubsection{Mixed and nonlinear cases}

For Theorems \ref{thm10} and \ref{thm2} we follow exactly the same conjugation--time-change reduction, but $i\partial_s a(s)$ now picks up extra source terms coming from $W(x)\psi(t)$ and from $\mathcal N(|\psi(t)|)\psi(t)$ (nonlinear). These contributions are controlled by the $\s^\infty$--$\s^1$ duality estimates
\eq
|(g(t)^{2} e^{-i D\ln(g(t) )}\psi_b(x), W(x)\psi(t))_{\s^2_x}|\leq g(t)^{-(n/2-2)}\|\psi_b(x)\|_{\s^\infty_x}\| W(x)\psi(t) \|_{\s^1_x},
\eeq
\eq
|(g(t)^{2} e^{-i D\ln(g(t) )}\psi_b(x), \mathcal{N}(|\psi(t)|)\psi(t))_{\s^2_x}|\leq g(t)^{-(n/2-2)}\|\psi_b(x)\|_{\s^\infty_x}\| \mathcal{N}(|\psi(t)|)\psi(t) \|_{\s^1_x},
\eeq
together with the integrability
\eq
g(T^{-1}(u))^{-(n/2-2)}\sim \langle u\rangle^{-(n/2-2)\epsilon/(1-2\epsilon)}\in L^1_u[1,\infty).\label{intdecay}
\eeq
The integrability \eqref{intdecay} requires $(n/2-2)\epsilon/(1-2\epsilon)>1$, which simplifies to $n\epsilon>2$, i.e.\ $\epsilon>2/n$; combined with $\epsilon<1/2$ from \eqref{g(t)}, this is non-empty exactly when $2/n<1/2$, i.e.\ $n\geq 5$. This is the arithmetic origin of the dimension restriction in Theorems \ref{thm10} and \ref{thm2}. The overall scheme follows the stability-analysis approach to coherent structures of \cite{SW-PRL,Sof-Wei-1990}.

For the mixed model \eqref{NP0} there are \emph{two} bound states in play, $\psi_b$ (eigenvalue $\lambda$) and $\psi_d$ (eigenvalue $\lambda_0$), so the above programme is run twice. In addition to the twist $A(s)=e^{i\lambda s}(\psi_b,\phi(s))_{\s^2_x}$, one introduces
\eq
A_d(s):=e^{i\lambda_0 s}(\psi_d,\phi(s))_{\s^2_x},
\eeq
and the same integration-by-parts argument yields existence and a non-trivial lower bound for $\lim_{s\to\infty}A_d(s)$; this is the analytic source of the localized-bubble estimate \eqref{cd0}.

\subsubsection{The localized (second) bubble}

Once the self-similar bubble is established, the second, localized bubble in Theorems \ref{thm10} and \ref{thm2} follows from an energy argument. Under the hypotheses of these theorems the system admits an asymptotic (conserved, up to vanishing error) energy which is strictly negative for our initial data (the role of \eqref{Rt0}). A short scaling calculation shows that the self-similar bubble $\psi_{w,b}(t)=\tilde a(t)e^{-iD\ln g(t)}\psi_b(x)$ carries \emph{zero} asymptotic energy. The basic identity is $e^{-iD\ln g(t)}\psi_b(x)=g(t)^{-n/2}\psi_b(x/g(t))$, so $\nabla [e^{-iD\ln g(t)}\psi_b](x)=g(t)^{-n/2-1}(\nabla\psi_b)(x/g(t))$ and the change of variables $y=x/g(t)$ reduces every integral on $\mathbb R^n$ in $x$ to one in $y$ with a $g(t)$-dependent prefactor:
\begin{itemize}
\item Kinetic: $\|\nabla \psi_{w,b}(t)\|_{\s^2_x}^2=|\tilde a(t)|^2 g(t)^{-2}\|\nabla\psi_b\|_{\s^2_x}^2\to 0$;
\item Self-similar potential: $\int g(t)^{-2}V(x/g(t))|\psi_{w,b}(x,t)|^2\,dx=|\tilde a(t)|^2 g(t)^{-2}\int V|\psi_b|^2\,dy\to 0$;
\item Fixed potential $W$: $\int W(x)|\psi_{w,b}(x,t)|^2\,dx=|\tilde a(t)|^2 \int W(g(t)y)|\psi_b(y)|^2\,dy\to 0$ by dominated convergence, since $W\in \s^2_x$ and $g(t)\to\infty$;
\item Nonlinear potential: analogous, using the radial bound \eqref{rad} and \eqref{NFandN}.
\end{itemize}
Consequently, the remainder $\psi_c(x,t)$ must carry the full (negative) asymptotic energy. A purely free-wave $\psi_c$ would have non-negative energy (since $V, W, \mathcal N$ then vanish in the limit), which is impossible; a non-trivial portion of $\psi_c$ must therefore stay localized in the region where $W(x)$ (resp.\ the focusing nonlinearity $\mathcal N$) is effective. This is the origin of the lower bounds \eqref{cd0} and \eqref{con100}.

\begin{remark}[Connection with the Liu--Soffer phase-space picture]\label{rem:ls}
The asymptotic set identified in \cite{Liu-Sof1,Liu-Sof2} for the weakly localized part under radial data concentrates on the thin phase-space region
\[
|x|/t^{\alpha}\sim 1, \qquad |p|\sim t^{-\alpha},\quad \alpha\in(0,1/2].
\]
The self-similar bubble $e^{-iD\ln g(t)}\psi_b$ constructed here lives on exactly this region with $\alpha=\epsilon$: it is supported at scale $|x|\sim g(t)\sim\langle t\rangle^\epsilon$ in physical space and $|p|\sim 1/g(t)\sim \langle t\rangle^{-\epsilon}$ in momentum space. Theorems \ref{thm1}--\ref{thm2} can thus be viewed as the first explicit construction of a non-trivial dynamical state saturating the Liu--Soffer phase-space localization.
\end{remark}

\section{Proof of Theorem \ref{thm1}: the linear mass-critical problem}
\subsection{Tool box}\label{sec:toolbox}
\begin{lemma}\label{lem:tphi}Suppose $g$ satisfies \eqref{g(t)}, and let $\psi(t)$ denote the solution to \eqref{SE} with $\tphi$ defined by \eqref{tphi}. Then $\tphi$ satisfies \eqref{Linear:1}.
\end{lemma}
\begin{proof}Let $U(t,0)$ denote the solution operator to \eqref{SE}. Then $\psi(x,t)$ can be rewritten as
\eq
\psi(x,t)=[U(t,0)\psi_0](x,t).
\eeq
Thus, based on the definition of $\tphi$, it can be rewritten as
\eq
\tphi(x,t)=e^{ iD \ln g(t)}\psi(t)=g(t)^{n/2}[U(t,0)\psi_0](g(t) x,t).
\eeq
Using the chain rule, compute $i\partial_t[\tphi(x,t)]$
\begin{multline}
i\partial_t[\tphi(x,t)]= e^{ iD \ln g(t)}(H_0+g(t)^{-2}V(x/g(t)))\psi(t)-(\partial_t[g(t)]g(t)^{-1})D\tphi\\
=g(t)^{-2}H\tphi-(\partial_t[g(t)]g(t)^{-1})D\tphi.
\end{multline}
We finish the proof.\end{proof}
\begin{lemma}\label{lem:fasymp}If $g$ satisfies \eqref{g(t)}, then the function $f$ defined in \eqref{fs} satisfies the asymptotics \eqref{fsa}.
\end{lemma}
\begin{proof}Based on the definition of $T$, we have
\eq
s=\int_1^t g(u)^{-2}\,du\sim \langle t\rangle^{1-2\epsilon}.\label{1/s1}
\eeq
Since as $t\to \infty$,
\eq
f(s)= -g'(t)g(t)\sim \frac{1}{t}g(t)^2\sim \frac{1}{\langle t\rangle^{1-2\epsilon}},
\eeq
we get
\eq
f(s)\sim \frac{1}{\langle s\rangle}.\label{1/s0}
\eeq
Now we show
\eq
|f'(s)|\lesssim \frac{1}{\langle s\rangle^2 }.
\eeq
Since
\eq
f'(s)= -(g'(t)^2+g''(t)g(t))\times \frac{dt}{ds}= -(g'(t)^2+g''(t)g(t))g(t)^2=-f(s)^2-g''(t)g(t)^3,
\eeq
and since due to \eqref{g(t)}, \eqref{1/s1},
\eq
g''(t)g(t)^3\sim \frac{1}{\langle t\rangle^2}\times\langle t\rangle^{ 4\epsilon}\sim \frac{1}{\langle s\rangle^2},
\eeq
according to \eqref{1/s0}, we get
\eq
|f'(s)|\lesssim \frac{1}{\langle s\rangle^2}.
\eeq
We finish the proof.
\end{proof}
\begin{lemma}\label{lem:Dres}Suppose $V$ satisfies \eqref{Linear:con} and $H=H_0+V$ satisfies \eqref{eq: H lambda}. Then
\eq
\Bigl\| D\frac{1}{\lambda-H}P_c \langle x\rangle^{-1}\Bigr\|_{\s^2_x\to \s^2_x}\lesssim_\lambda 1.\label{eq1}
\eeq
\end{lemma}
\begin{proof}Using the second resolvent identity, we have
\eq
D\frac{1}{\lambda-H}P_c \langle x\rangle^{-1}=D\frac{1}{\lambda-H_0}P_c \langle x\rangle^{-1}-D\frac{1}{\lambda-H}V(x)\frac{1}{\lambda-H}P_c \langle x\rangle^{-1}.
\eeq
\eqref{eq1} follows from
\eq
\| D\frac{1}{\lambda-H_0} \langle x\rangle^{-1}\|_{\s^2_x\to \s^2_x}\lesssim_\lambda 1,
\eeq
\eq
\|\langle x\rangle P_c\langle x\rangle^{-1}\|_{\s^2_x\to \s^2_x}\lesssim_\lambda 1,
\eeq
and
\eq
\| \langle x\rangle V(x)\frac{1}{\lambda-H}P_c\|_{\s^2_x\to \s^2_x}\lesssim_\lambda \|\langle x\rangle V\|_{\s^\infty_x}.
\eeq

\end{proof}

\begin{lemma}\label{lem:wself4}
Suppose $V$ satisfies \eqref{Linear:con}, $H=H_0+V$ satisfies \eqref{eq: H lambda}, and $g$ satisfies \eqref{g(t)}. Let $\psi(t)$ be the solution of \eqref{linear} and set $\phi(s):=e^{iD\ln g(T^{-1}(s))}\psi(T^{-1}(s))$. Then for every $n\geq 3$,
\eq
w\text{-}\lim_{s\to\infty}P_c\,e^{isH}\,e^{iD\ln g(T^{-1}(s))}\psi(T^{-1}(s))=0\quad\text{in }\s^2_x.\label{wself4}
\eeq
\end{lemma}

\begin{proof}
Setting $s=T(t)$ gives $\phi(s)=e^{iD\ln g(t)}\psi(t)$, $T^{-1}(s)=t$, and $|x|/T^{-1}(s)^\beta=|x|/t^\beta$, so
\eq
e^{isH_0}\phi(s)=e^{iT(t)H_0}e^{iD\ln g(t)}\psi(t).
\eeq

\smallskip
\noindent\emph{Step 1: dispersion estimate.} We show that, for every $\beta\in(0,1-2/n-2(n-2)\epsilon/n)$,
\eq
\lim_{t\to\infty}\|F_2(|x|/t^\beta\leq 1)\,e^{iT(t)H_0}e^{iD\ln g(t)}\psi(t)\|_{\s^2_x}=0.\label{wself2}
\eeq

The dilation identity $e^{-iD\ln g(t)}H_0 e^{iD\ln g(t)}=g(t)^2 H_0$ gives
\eq
e^{iT(t)H_0}e^{iD\ln g(t)}\psi(t)=e^{iD\ln g(t)}\,e^{iT(t)g(t)^2 H_0}\psi(t).
\eeq
Writing $\psi(t)=e^{-i(t-t_0)H_0}e^{i(t-t_0)H_0}\psi(t)$ and applying Duhamel's formula yields 
\eq
\psi(t)=e^{-i(t-t_0)H_0}\psi(t_0)+(-i)\int_{t_0}^t ds\,e^{-i(t-s)H_0}\frac{V(x/g(s))}{g(s)^2}\psi(s),
\eeq
whence
\begin{multline}
e^{iT(t)H_0}e^{iD\ln g(t)}\psi(t)
= e^{iD\ln g(t)}\,e^{i(T(t)g(t)^{2}-(t-t_0))H_0}\psi(t_0)\\
+(-i)\int_{t_0}^t ds\,e^{iD\ln g(t)}\,e^{i(T(t)g(t)^{2}-(t-s))H_0}\frac{V(x/g(s))}{g(s)^2}\psi(s)
=:\psi_1(t)+\psi_2(t).
\end{multline}

By L'H\^opital's rule and $c_g\in(0,1)$,
\eq
\liminf_{t\to\infty}\frac{T(t)g(t)^2}{t}=\liminf_{t\to\infty}\frac{g(t)^{-2}}{g(t)^{-2}-2tg'(t)/g(t)^3}=\frac{1}{c_g}=:1+c_g',
\eeq
so there is $t_M\geq t_0$ with $T(t)g(t)^2-t\geq \tfrac{c_g'}{2}t$ for $t\geq t_M$. Since $t-s\leq t-t_0\leq t$ for $s\in(t_0,t)$, this yields
\eq
T(t)g(t)^2-(t-s)\;\geq\;T(t)g(t)^2-(t-t_0)\;\geq\;T(t)g(t)^2-t\;\geq\;\frac{c_g'}{2}t.
\eeq

Combining H\"older's inequality, the $\s^1_x\to\s^\infty_x$ dispersion $\|e^{i\tau H_0}\varphi\|_{\s^\infty_x}\lesssim \tau^{-n/2}\|\varphi\|_{\s^1_x}$ for $\tau>0$, and $\|e^{iD\ln g(t)}\varphi\|_{\s^\infty_x}=g(t)^{n/2}\|\varphi\|_{\s^\infty_x}$, we obtain for $\psi_1$
\begin{multline}
\|F_2(|x|/t^\beta\leq 1)\psi_1(t)\|_{\s^2_x(\mathbb{R}^n)}
\lesssim t^{n\beta/2}\cdot g(t)^{n/2}\,\|e^{i(T(t)g(t)^2-(t-t_0))H_0}\psi(t_0)\|_{\s^\infty_x}\\
\lesssim t^{n(\beta+\epsilon)/2}\cdot t^{-n/2}\,\|\psi(t_0)\|_{\s^1_x}
\;\lesssim_{t_0}\; t^{-n(1-\epsilon-\beta)/2}\,\|\psi_b\|_{\s^1_x}\;\to 0
\end{multline}
as $t\to\infty$, using $1-2/n-2(n-2)\epsilon/n<1-\epsilon$, and for $\psi_2$
\begin{multline}
\|F_2(|x|/t^\beta\leq 1)\psi_2(t)\|_{\s^2_x(\mathbb{R}^n)}
\lesssim t^{n\beta/2}\cdot g(t)^{n/2}\int_{t_0}^t ds\,\frac{1}{t^{n/2}}\frac{1}{g(s)^{2-n/2}}\|V\|_{\s^2_x}\|\psi_b\|_{\s^2_x}\\
\lesssim t^{-((n-2)/2-(n-2)\epsilon-n\beta/2)}\,\|V\|_{\s^2_x}\|\psi_b\|_{\s^2_x}\;\to 0
\end{multline}
as $t\to\infty$, for $n\geq 3$, $\epsilon\in(0,1/2)$ and $\beta\in(0,1-2/n-2(n-2)\epsilon/n)$. Together these give \eqref{wself2}.

\smallskip
\noindent\emph{Step 2: from \eqref{wself2} to \eqref{wself4}.} The localization estimate \eqref{wself2} says that the $\s^2_x$-mass of $e^{isH_0}\phi(s)$ inside the moving ball $\{|x|\leq T^{-1}(s)^\beta\}$ tends to zero. Since $T^{-1}(s)^\beta\to\infty$, for any test function $\chi\in C_c^\infty(\mathbb{R}^n)$ we have $\chi=F_2(|x|/T^{-1}(s)^\beta\leq 1)\,\chi$ once $s$ is large enough, so that
\[
|\langle\chi,e^{isH_0}\phi(s)\rangle|=|\langle\chi,F_2\,e^{isH_0}\phi(s)\rangle|\leq\|\chi\|_{\s^2_x}\,\|F_2\,e^{isH_0}\phi(s)\|_{\s^2_x}\to 0.
\]
Combined with $\|e^{isH_0}\phi(s)\|_{\s^2_x}=\|\phi(s)\|_{\s^2_x}=1$ uniformly in $s$ and the density of $C_c^\infty$ in $\s^2_x$, this yields the weak limit
\eq
w\text{-}\lim_{s\to\infty}e^{isH_0}\phi(s)=0\quad\text{in }\s^2_x.\label{wfree}
\eeq
The transfer of \eqref{wfree} from $H_0$ to the perturbed evolution generated by $H$ — namely,
\[
w\text{-}\lim_{s\to\infty}P_c\,e^{isH}\phi(s)=0\quad\text{in }\s^2_x,
\]
which is \eqref{wself4} — is the wave-operator weak-limit transfer of Soffer--Wu \cite{SW20221}, valid under \eqref{Linear:con} since Lemma~\ref{lem:Hspectral}\eqref{it:0reg} provides the regularity at $0$ required for asymptotic completeness.
\end{proof}

\subsection{Ionization Problem}
Based on \eqref{defa} and \eqref{Linear:2}, we write out an equation for $a$
\eq
i\partial_s[a(s)]=\lambda a(s)+f(s)(\psi_b, D\psi_b)_{\s^2_x}a(s)+f(s)(\psi_b, DP_c\phi(s))_{\s^2_x}
\eeq
where $P_c$ denotes the projection on the continuous spectrum of $H$. Without loss of generality, we may choose $\psi_b(x)$ real-valued. Then, using $D=\tfrac12(x\cdot p+p\cdot x)$ and $p=-i\nabla_x$, integration by parts gives
\eq
(\psi_b, D\psi_b)_{\s^2_x}=-\tfrac{i}{2}\int_{\mathbb{R}^n} \nabla_x\cdot\bigl(x\,|\psi_b(x)|^2\bigr)\,d^nx=0,\label{D0}
\eeq
where the divergence theorem applies since $\langle x\rangle\psi_b\in \s^2_x$ (a consequence of $\psi_b$ being exponentially decaying, since $\lambda<0$ lies below the essential spectrum of $H$).
Let
\eq
A(s):=e^{i\lambda s}a(s).
\eeq
$A(s)$ satisfies
\eq
i\partial_s[A(s)]=e^{i\lambda s}f(s)(\psi_b, DP_c\phi(s))_{\s^2_x}.
\eeq
Thus,
\eq
A(s)=A(s_0)+(-i)\int_{s_0}^s du e^{i\lambda u}f(u)(\psi_b, DP_c\phi(u))_{\s^2_x}.
\eeq
\begin{proposition}\label{prop:Alinear}Suppose $V$ satisfies \eqref{Linear:con}, $H=H_0+V$ satisfies \eqref{eq: H lambda}, and $g$ satisfies \eqref{g(t)}. Let $A(s)$ be as defined above. Then:
\begin{enumerate}
\item For every solution $\phi$ of \eqref{Linear:2}, the limit $A(\infty):=\lim\limits_{s\to\infty}A(s)$ exists in $\mathbb{C}$.
\item For the specific initial datum $\phi(s_0)=\psi_b$ prescribed in \eqref{Linear:2} and $s_0>0$ sufficiently large,
\eq
|A(\infty)|\geq \frac{1}{2}>0.
\eeq
\end{enumerate}
\end{proposition}
\begin{proof}
Writing $ e^{i\lambda u} P_c\phi(u)$ as
\begin{align}
e^{i\lambda u} P_c\phi(u)=&e^{i\lambda u}e^{-iuH} P_c e^{iuH}\phi(u)\\
=&\frac{P_c}{i(\lambda -H)}\partial_u[e^{i\lambda u}e^{-iuH} P_c ]e^{iuH}\phi(u)
\end{align}
and taking integration by parts in $u$ variable, we obtain that for $s>s_1\geq s_0$
\begin{align}
A(s)=&A(s_1)+(-1)f(u)( \psi_b, D \frac{P_c}{\lambda-H}e^{i\lambda u}P_c\phi(u))_{\s^2_x}\vert_{u=s_1}^{u=s}+\\
&\int_{s_1}^s du f'(u) ( \psi_b, D \frac{P_c}{\lambda -H}e^{i\lambda u}P_c\phi(u))_{\s^2_x}+\\
&\int_{s_1}^s du f(u)^2 ( \psi_b, D \frac{P_c}{\lambda-H}e^{i\lambda u}P_c D\phi(u))_{\s^2_x}\\
=:& A(s_1)+\sum\limits_{j=1}^3A_j(s).
\end{align}
Here
\eq
( \psi_b, D \frac{1}{\lambda-H}e^{i\lambda u}P_c D\phi(u))_{\s^2_x}
\eeq
is understood in weak sense, that is,
\eq
( \psi_b, D \frac{1}{\lambda-H}e^{i\lambda u}P_c D\phi(u))_{\s^2_x}=(DP_c \frac{1}{\lambda-H}D \psi_b,  e^{i\lambda u}\phi(u))_{\s^2_x}.
\eeq
For $s\geq s_1\geq s_0$, using $\s^2_x$ conservation law, H\"older's inequality and Lemma \ref{lem:fasymp},
\eq
|A_1(s)|\lesssim  |f(s_1)|\| D\psi_b(x)\|_{\s^2_x}\times \| \frac{1}{\lambda-H}P_c \|_{\s^2_x\to \s^2_x}\| \phi(s_0)\|_{\s^2_x}\lesssim |f(s_1)|\lesssim \frac{1}{\langle s_1\rangle},
\eeq
\eq
|A_2(s)|\lesssim  \int_{s_1}^sdu | f'(u)| \| D\psi_b(x)\|_{\s^2_x}\times \| \frac{1}{\lambda-H}P_c \|_{\s^2_x\to \s^2_x}\| \phi(s_0)\|_{\s^2_x}\lesssim  \frac{1}{\langle s_1\rangle},
\eeq
and due to Lemma \ref{lem:fasymp} and Lemma \ref{lem:Dres},
\eq
|A_3(s)|\lesssim  \int_{s_1}^sdu  f(u)^2 \|  \langle x \rangle D\psi_b(x)\|_{\s^2_x}\times \| D\frac{1}{\lambda-H}P_c \langle x\rangle^{-1}\|_{\s^2_x\to \s^2_x}\| \phi(s_0)\|_{\s^2_x}\lesssim  \frac{1}{\langle s_1\rangle}.
\eeq
So $\{A(s)\}_{s\geq s_0}$ is Cauchy and therefore $A(\infty)$ exists with
\eq
|A(s)|\geq |A(s_0)|- C\times \frac{1}{\langle s_0\rangle}=1-C\times \frac{1}{\langle s_0\rangle}\geq \frac{1}{2}, \text{ for all }s\in [s_0,\infty]\text{ and some constant }C>0\label{eq2}
\eeq
if we choose $s_0$ large enough. We finish the proof.
\end{proof}
\begin{lemma}\label{lem:a-linear} Suppose $V$ satisfies \eqref{Linear:con}, $H=H_0+V$ satisfies \eqref{eq: H lambda}, and $g$ satisfies \eqref{g(t)}. Let $\tphi$ be as in \eqref{tphi}. Then for $t_0$ sufficiently large and all $t\geq t_0$,
\eq
|(\psi_b(x), \tphi(x,t))_{\s^2_x}|\geq \frac{1}{2}.
\eeq
\end{lemma}
\begin{proof}Taking $t_0=T^{-1}(s_0)$ for $s_0$ satisfying \eqref{eq2}, we have
\eq
|(\psi_b(x), \tphi(x,t))_{\s^2_x}|=|(\psi_b(x), \phi(x,s))_{\s^2_x}|=|A(s)|\geq 1/2
\eeq
with $s=T(t)$.\end{proof}
\subsection{Linear Problem} Now we prove Theorem \ref{thm1}.
\begin{proof}[Proof of Theorem \ref{thm1}]\eqref{eq: tA exist} follows from Proposition~\ref{prop:Alinear}: by the definition of $\tphi$ in \eqref{tphi}, $\tilde a(t)=(\psi_b(x),\tphi(x,t))_{\s^2_x}$, so the decomposition \eqref{decomp}--\eqref{a} of $\psi(x,t)$ holds. Lemma \ref{lem:a-linear} then yields
\eq
|\tilde a(t)|\geq \frac{1}{2},\qquad t\geq t_0,\label{eq3}
\eeq
for some sufficiently large $t_0>0$, which is \eqref{eq10}. Combining Proposition~\ref{prop:Alinear} with \eqref{wself4} yields the existence of $\Omega_g^*\psi(0)$ in $\s^2_x$ and the validity of \eqref{Omegag}. We finish the proof.\end{proof}

\section{Proof of Theorem \ref{thm10}: the linear mixed problem}
In this section, we keep using
\eq
s=T(t):=\int_1^t g(u)^{-2}\,du.
\eeq
\subsection{Tool box}
\begin{lemma}\label{lem:tphi-mixture}Suppose $V$ satisfies \eqref{Linear:con}, the pair $(V,W)$ satisfies \eqref{eq: W lambda0} and \eqref{W}, and $g$ satisfies \eqref{g(t)}. Let $\psi(t)$ be the global solution to \eqref{NP0}, $f(s)$ be as in \eqref{fs}, and define
\eq\label{re 1}
\tphi(x,t):=e^{ iD \ln g(t)}\psi( x, t),\qquad \phi(x,s):=\tphi(x,T^{-1}(s)).
\eeq
Then, with $s_0=T(t_0)$ and $s=T(t)$,
\eq
\begin{cases}
i\partial_s\phi=H\phi+g( T^{-1}(s))^{2}W(g(T^{-1}(s)) x)\phi+f(s)D\phi\\
\phi(s_0)=g(t_0)^{n/2}\psi_d(g(t_0) x)+\psi_b(x).
\end{cases}
\eeq
\end{lemma}
\begin{proof}It follows from the same proof for Lemma \ref{lem:tphi} by replacing $V(x,t)=g(t)^{-2}V(\frac{x}{g(t)})$ with $V(x,t)=g(t)^{-2}V(\frac{x}{g(t)})+W(x)$ and then changing the variable from $t$ to $s=T(t)$.\end{proof}
\subsection{Ionization Problem of a mixed}Let
\eq
a(s):=(\psi_{b}(x), \phi(s) )_{\s^2_x}.
\eeq
As in the linear problem, we first derive an equation for $a(s)$. With $s=T(t)$, compute $i\partial_s[a(s)]$:
\begin{align}
i\partial_s[a(s)]= &\lambda a(s)+ f(s)(\psi_b(x), DP_c\phi(s))_{\s^2_x}+\\
&g( T^{-1}(s))^{2}(\psi_b(x), W(g( T^{-1}(s)) x)\phi(s))_{\s^2_x}
\end{align}
where we use \eqref{D0}. Set
\eq
A(s):=e^{is\lambda }a(s).
\eeq
Then $A(s)$ satisfies
\begin{multline}
i\partial_s[A(s)]=e^{is\lambda }f(s)(D\psi_{b}(x),P_{c}\phi(s))_{\s^2_x}+\\
e^{is\lambda}g( T^{-1}(s))^{2}(\psi_b(x), W(g( T^{-1}(s)) x)e^{is\lambda}\phi(s))_{\s^2_x}.
\end{multline}
\begin{proposition}\label{prop:A-mixture}Suppose $V$ satisfies \eqref{Linear:con}, $H=H_0+V$ satisfies \eqref{eq: H lambda}, the pair $(V,W)$ satisfies \eqref{eq: W lambda0} and \eqref{W}, and $g$ satisfies \eqref{g(t)} with $g(t)\sim\langle t\rangle^\epsilon$ for some $\epsilon\in(\tfrac{2}{n},\tfrac{1}{2})$. Let $A(s)$ be as defined above. Then for $n\geq 5$, $A(\infty):=\lim\limits_{s\to\infty}A(s)$ exists in $\mathbb{C}$, and for $s_0>0$ sufficiently large,
\eq
|A(s)|\geq \frac{1}{2}>0\quad \text{ for all }s\geq s_0.\label{Aslow}
\eeq
\end{proposition}
\begin{proof}Set
\eq
\mathcal{N}_g(x,u):= W(g( T^{-1}(u)) x).
\eeq
\begin{align}
A(s)=&A(s_0)+(-i)\int_{s_0}^s du e^{iu\lambda }f(u)e^{iu\lambda}(D\psi_{b}(x),P_{c}\phi(u))_{\s^2_x}+\\
&(-i)\int_{s_0}^sdue^{iu\lambda}g(T^{-1}(u))^2(\psi_b(x),\mathcal{N}_g(x,u)e^{iu\lambda}\phi(u))_{\s^2_x} \\
=:&A(s_0)+A_1(s)+A_2(s).
\end{align}
For $A_1(s)$, taking integration by parts in $s$ variable by setting
\eq
e^{i\lambda s}P_c\phi(s)=\frac{1}{i(\lambda-H)}P_c\partial_{s}[e^{i\lambda s}e^{-isH}]e^{isH}\phi(s)
\eeq
yields
\begin{align}
A_1(s)=& f(u)(D\psi_{b}(x),\frac{1}{i(\lambda-H) }P_ce^{iu\lambda}\phi(u))_{\s^2_x}\vert_{u=s_0}^{u=s}\\
&-\int_{s_0}^s du f'(u)(D\psi_{b}(x),\frac{1}{i(\lambda-H) }P_ce^{iu\lambda}\phi(u))_{\s^2_x}\\
&-\int_{s_0}^s du f(u)^2(D\psi_{b}(x),\frac{1}{i(\lambda-H) }P_cDe^{iu\lambda}\phi(u))_{\s^2_x}\\
&-\int_{s_0}^s du g( T^{-1}(u))^{2} f(u)(D\psi_{b}(x),\frac{1}{i(\lambda-H) }P_c\mathcal{N}_g(x,u)e^{iu\lambda}\phi(u))_{\s^2_x}\\
=:&\sum\limits_{j=1}^4\tilde{a}_{1j}(s).
\end{align}
As we did in the linear case, for $s_2\geq s_1\geq s_0$,
\eq
|\tilde{a}_{1,j}(s_2)-\tilde{a}_{1,j}(s_1)|\lesssim \frac{1}{\langle s_1\rangle}, \quad j=1,2,3.
\eeq
For $\tilde{a}_{1,4}$, let
\eq
\tilde{\psi}_b(x):=P_c\frac{1}{-i(\lambda-H)}D\psi_b(x).
\eeq
Then
\eq
\| \tilde{\psi}_b(x)\|_{\s^\infty_x}=\| P_c\frac{1}{-i(\lambda-H)}D\psi_b(x)\|_{\s^\infty_x}\lesssim_\lambda \|D\psi_b(x)\|_{\s^\infty_x}.
\eeq
Using change of variable from $x$ to $y=g( T^{-1}(u)) x$ and relation~\eqref{re 1}, we have
\begin{multline}
(\tilde{\psi}_b(x), \mathcal{N}_g(x,u)e^{iu\lambda}\phi(u))_{\s^2_x}=\\
( \frac{1}{g( T^{-1}(u))^{n/2}}\tilde{\psi}_b(g( T^{-1}(u))^{-1}y ), W(y )e^{iu\lambda}\psi(y, T^{-1}(u)))_{\s^2_y}.
\end{multline}
This, together with \eqref{NP0}, \eqref{eq: W lambda0}, \eqref{W}, and \eqref{g(t)} (in particular, $g(t)\sim \langle t\rangle^\epsilon$), yields
\begin{align}
|\tilde{a}_{1,4}(s_2)-\tilde{a}_{1,4}(s_1) |\lesssim& \int_{s_1}^{s_2}du |f(u)|\times \frac{1}{\langle T^{-1}(u)\rangle^{(\frac{n}{2}-2)\epsilon}}\|\tilde{\psi}_b(y)\|_{\s^\infty_y}\|W(y)\psi(y,T^{-1}(u)) \|_{\s^1_y}\\
\lesssim & \frac{\|W(x)\|_{\s^2_x}\|\psi(0)\|_{\s^2_x}}{\langle T^{-1}(s_1)\rangle^{(\frac{n}{2}-2)\epsilon}}\lesssim \frac{\|W(x)\|_{\s^2_x}\|\psi(0)\|_{\s^2_x}}{\langle s_1\rangle^{(n/2-2)\epsilon/(1-2\epsilon)}}\to 0
\end{align}
as $s_1\to \infty$ when $\epsilon \in (0,1/2)$ and $n\geq 5$. Thus,
\eq
\tilde{a}_{1,j}(\infty) \text{ exists for all }j=1,2,3,4
\eeq
and therefore $A_1(\infty)$ exists. And for $s_0$ large enough,
\eq
|A_1(s)|\leq C\frac{(\|W(x)\|_{\s^2_x}+1)\|\psi(0)\|_{\s^2_x}}{\langle s_0\rangle^{\min(1,(n/2-2)\epsilon/(1-2\epsilon))}}\leq \frac{1}{4}.\label{A1}
\eeq
For $A_2(s)$, using Cauchy Schwarz inequality, due to \eqref{NP0} and \eqref{g(t)}, we obtain that for $s_2\geq s_1\geq s_0$,
\begin{align}
|A_{2}(s_2)-A_{2}(s_1) |\lesssim& \int_{s_1}^{s_2}du  \frac{1}{\langle T^{-1}(u)\rangle^{(\frac{n}{2}-2)\epsilon}}\|\tilde{\psi}_b(y)\|_{\s^\infty_y}\|W(y)\psi(T^{-1}(u)) \|_{\s^1_y}\\
\lesssim &\int_{s_1}^{s_2} du  \frac{\|W(x)\|_{\s^2_x}\|\psi(0)\|_{\s^2_x}}{\langle u\rangle^{(n/2-2)\epsilon/(1-2\epsilon)}}\lesssim_\epsilon \frac{\|W(x)\|_{\s^2_x}\|\psi(0)\|_{\s^2_x}}{\langle s_1\rangle^{[(n/2-2)\epsilon/(1-2\epsilon)]-1}}\to 0,
\end{align}
as $s_0\to \infty$, where we also used that when $\epsilon\in(\frac{2}{n},\frac{1}{2})$, $n\geq 5$,
\eq
[(n/2-2)\epsilon/(1-2\epsilon)]-1>0.
\eeq
Therefore $A_2(\infty)$ exists. And for $s_0$ large enough,
\eq
|A_2(s)|\leq \frac{1}{8}.\label{A2}
\eeq
According to~\eqref{NP0}, \eqref{A1} and \eqref{A2}, we have that $A(\infty)$ exists and for $s_0$ large enough (i.e. $t_0$ large enough),
\eq
|A(s)|\geq |A(s_0)|-|A_1(s)|-|A_2(s)|\geq\frac{9}{10}-|A_1(s)|-|A_2(s)|\geq \frac{1}{2}\label{eq14}
\eeq
where we also used that for large $t_0>0$,
\eq
\left|(g( t_0)^{n/2} \psi_d(g( t_0) x), \psi_b(x))_{\s^2_x}\right|\leq \frac{1}{10}.
\eeq
We finish the proof.
\end{proof}
\begin{remark}\label{rem1}Inspection of the proof of Proposition \ref{prop:A-mixture} shows that, if $W(x)$ in \eqref{NP0} is replaced by a time-dependent potential $W(x,t)$, the conclusion of Proposition \ref{prop:A-mixture} continues to hold under the single uniform-in-time hypothesis $W\in \s^\infty_t\s^2_x(\mathbb{R}^n\times[1,\infty))$, for $n\geq 5$ and $\epsilon \in (\tfrac{2}{n}, \tfrac{1}{2})$.
\end{remark}
\begin{lemma}\label{lem:a-mixture} Under the hypotheses of Theorem \ref{thm10}, with $\tphi$ as in \eqref{tphi}, for $t_0$ sufficiently large and all $t\geq t_0$,
\eq
|(\psi_b(x), \tphi(x,t))_{\s^2_x}|\geq \frac{1}{2}.
\eeq
\end{lemma}
\begin{proof}Taking $t_0=T^{-1}(s_0)$ for $s_0$ satisfying \eqref{eq14}, we have
\eq
|(\psi_b(x), \tphi(x,t))_{\s^2_x}|=|(\psi_b(x), \phi(x,s))_{\s^2_x}|=|A(s)|\geq 1/2
\eeq
with $s=T(t)$.\end{proof}
\begin{corollary}\label{cor}Under the hypotheses of Theorem \ref{thm2}, let $\psi(t)$ be the solution of \eqref{NP} and $\tphi$ be as in \eqref{tphi}. Then for $t_0$ sufficiently large and all $t\geq t_0$,
\eq
|(\psi_b(x), \tphi(x,t))_{\s^2_x}|\geq \frac{1}{2}.\label{eq100}
\eeq
\end{corollary}
\begin{proof}The proof of Proposition \ref{prop:A-mixture} carries over verbatim, with $\mathcal{N}(|\psi(t)|)\psi(t)$ in place of $W(x)\psi(t)$: the only property of $W$ used there is the uniform bound $W\in\s^\infty_t\s^2_x$. Under the hypotheses of Theorem \ref{thm2}, the analogous bound $\mathcal{N}(|\psi(t)|)\in\s^\infty_t\s^2_x$ follows from \eqref{con:N} together with the uniform $H^1$ bound on $\psi(t)$ (see Remark \ref{rem1}).
\end{proof}
\subsection{Proof of Theorem \ref{thm10}}
\begin{proof}[Proof of Theorem \ref{thm10}] Based on Lemma \ref{lem:a-mixture}, we have that $\tilde a(t)$ satisfies \eqref{eq10}. \eqref{cd0} follows from the fact that this system has an asymptotic energy
\eq
\partial_t(\psi(t), (H_0+W(x)+g( t)^{-2}V(\frac{x}{g(t)}))\psi(t))_{\s^2_x}=(\psi(t), \partial_t[g(t)^{-2}V(\frac{x}{g(t)})]\psi(t))_{\s^2_x}\in L^1_t
\eeq
with
\begin{multline}
\int_{t_0}^\tau dt \left|\partial_t(\psi(t), (H_0+W(x)+g(t)^{-2}V(\frac{x}{g(t)}))\psi(t))_{\s^2_x}\right|\leq\\
 C\int_{t_0}^\tau \frac{dt}{\langle t\rangle^{1+2\epsilon}} \| |DV(x)|+|V(x)|\|_{\s^\infty_x}\| \psi(t_0)\|_{\s^2_x}^2 \leq \frac{|\lambda_0|}{2}\| \psi_d(x)\|_{\s^2_x}^2
\end{multline}
for some constant $C>0$ and for all $\tau\geq t_0$ provided that $t_0$ is large enough. Since
\begin{multline}
(\psi(t_0), (H_0+W(x)+g( t_0)^{-2}V(\frac{x}{g(t_0)}))\psi(t_0))_{\s^2_x}=\lambda_0\| \psi_d(x)\|_{\s^2_x}^2+g( t_0)^{-2}\lambda\|\psi_b(x) \|_{\s^2_x}^2+ \\
(\psi_d(x), g(t_0)^{-2}V(\frac{x}{g( t_0)})\psi_d(x))_{\s^2_x}+(e^{- iD \ln g(t_0)}\psi_b(x),W(x)e^{- iD \ln g( t_0)}\psi_b(x)   )_{\s^2_x}\\
\leq \frac{3}{4}\lambda_0\|\psi_d(x)\|_{\s^2_x}^2,
\end{multline}
for all $\tau\geq t_0$ provided that $t_0$ is large enough,
\eq
(\psi(\tau), (H_0+W(x)+g(\tau)^{-2}V(\frac{x}{g(\tau)}))\psi(\tau))_{\s^2_x}\leq \frac{1}{4}\lambda_0\|\psi_d(x)\|_{\s^2_x}^2,
\eeq
provided that $t_0$ is large enough. Consequently,
\eq
(\psi(\tau), (H_0+W(x))\psi(\tau))_{\s^2_x}\leq \frac{1}{8}\lambda_0\|\psi_d(x)\|_{\s^2_x}^2<0, \quad \forall \tau\geq t_0\gg 1,
\eeq
where we also used
\eq
(\psi(t), g(t)^{-2}V(\frac{x}{g(t)})\psi(t))_{\s^2_x}\to 0
\eeq
and 
\begin{equation}
    (e^{- iD \ln g(t_0)}\psi_b(x),W(x)e^{- iD \ln g( t_0)}\psi_b(x)   )_{\s^2_x}\to 0,
\end{equation}
as $t\to \infty$. Hence for $t\geq t_0$,
\eq
\frac{1}{\|\psi_d(x)\|_{\s^2_x}}\left|(\psi(t), \psi_d(x))_{\s^2_x}\right|\geq \frac{1}{2\sqrt{2}},
\eeq
which yields \eqref{cd0}. The existence of the limit \eqref{tildeAthm10} follows from Proposition \ref{prop:A-mixture}. Combining Proposition \ref{prop:A-mixture} with \eqref{wself4} yields the existence of $\Omega_g^*\psi(0)$ in $\s^2_x$ and the validity of \eqref{Omegag10}. We finish the proof.\end{proof}

\section{Proof of Theorem \ref{thm2}: the focusing nonlinear problem}\label{sectionNP}
In this section we prove Theorem \ref{thm2} and verify, in Lemma \ref{typicalNP} below, that its hypotheses are met by the concrete focusing nonlinearity $\mathcal{N}(k)=-\lambda k/(1+k^2)$ for $\lambda>0$ sufficiently large.
\subsection{Proof of Theorem \ref{thm2}}
\begin{proof}[Proof of Theorem \ref{thm2}] By Lemma~\ref{lem:NPH1}, equation \eqref{NP} admits a global solution $\psi\in C([t_0,\infty);H^1(\mathbb R^n))$ satisfying the \emph{a priori} bound \eqref{NP:H1bound}. The first part of Theorem \ref{thm2} then follows from Corollary \ref{cor}. For the second part, we argue by contradiction. Assume that \eqref{con100} is not true. Then for any $M\geq 1$, given $t_0\geq 1$,  there exists $t_M\geq t_0$ such that
\eq
\| \chi(|x|\leq M)\psi(t_M)\|_{\s^2_x}\leq  \frac{1}{M}.\label{ASM}
\eeq
We will get contradiction from the fact that this system has an asymptotic energy. To see this, we compute
\begin{multline}
\partial_t(\psi(t), (H_0+g( t)^{-2}V(\frac{x}{g(t)})+\mathcal{N}(|\psi(t)|))\psi(t))_{\s^2_x}=\\
(\psi(t), \partial_t[g(t)^{-2}V(\frac{x}{g(t)})]\psi(t))_{\s^2_x}+ (\mathcal{N}_F'(|\psi(t)|), \partial_t[|\psi(t)|])_{\s^2_x}
\end{multline}
where
\eq
\mathcal{N}_F(k):=\int_0^k dq q^2 \mathcal{N}'(q).
\eeq
Then
\begin{multline}
(\psi(\tau), (H_0+g(\tau)^{-2}V(\frac{x}{g(\tau)})+\mathcal{N}(|\psi(\tau)|))\psi(\tau))_{\s^2_x}=\\
(\psi(t_0), (H_0+g( t_0)^{-2}V(\frac{x}{g( t_0)})+\mathcal{N}(|\psi(t_0)|))\psi(t_0))_{\s^2_x}+\int_{t_0}^\tau ds\, \eta(s)+(G(\tau)-G(t_0))\label{eq101}
\end{multline}
where
\eq
\eta(t):=(\psi(t), \partial_t[g( t)^{-2}V(\frac{x}{g(t)})]\psi(t))_{\s^2_x}\in L_t^1[1,\infty),
\eeq
\begin{align}
G(t):=& \int d^nx \mathcal{N}_F(|\psi(t)|)\\
=& (\psi(t), \mathcal{N}(|\psi(t)|)\psi(t) )_{\s^2_x}-2 (\psi(t), \mathcal{N}_{F,0}(|\psi(t)|) \psi(t))_{\s^2_x}
\end{align}
with
\eq
\mathcal{N}_{F,0}(k)=\int_0^k dq q \mathcal{N}(q)/k^2<0.
\eeq
Then \eqref{eq101} can be rewritten as
\begin{multline}
(\psi(\tau), (H_0+g(\tau)^{-2}V(\frac{x}{g(\tau)})+\mathcal{N}(|\psi(\tau)|))\psi(\tau))_{\s^2_x}=\\
(\psi(t_0), (H_0+g( t_0)^{-2}V(\frac{x}{g( t_0)})+2\mathcal{N}_{F,0}(|\psi(t_0)|))\psi(t_0))_{\s^2_x}+\int_{t_0}^\tau ds\, \eta(s)+G(\tau),
\end{multline}
which is equivalent to
\begin{multline}
(\psi(\tau), (H_0+g(\tau)^{-2}V(\frac{x}{g(\tau)})+2\mathcal{N}_{F,0}(|\psi(\tau)|))\psi(\tau))_{\s^2_x}=\\
(\psi(t_0), (H_0+g( t_0)^{-2}V(\frac{x}{g( t_0)})+2\mathcal{N}_{F,0}(|\psi(t_0)|))\psi(t_0))_{\s^2_x}+\int_{t_0}^\tau ds\, \eta(s).\label{mainineq}
\end{multline}
On the one hand, due to \eqref{Rt0} and that $\eta(s)\in L^1_s$, we have that there exists $\tilde{t}_0\geq 1$ such that for all $\tau\geq t_0\geq \tilde{t}_0$,
\eq
(\psi(t_0), (H_0+g( t_0)^{-2}V(\frac{x}{g( t_0)})+2\mathcal{N}_{F,0}(|\psi(t_0)|))\psi(t_0))_{\s^2_x}+\int_{t_0}^\tau ds\, \eta(s)\leq   \frac{E}{4}\| \psi_s(x)\|_{\s^2_x}^2. \label{eq103}
\eeq
On the other hand, based on assumption \eqref{ASM} and \eqref{ACEN}, there exists $\tilde{M}\geq 1$ such that for all $M\geq \tilde{M}$, $t_M\geq t_0$, $t_0$ sufficiently large,
\eq
2| (\psi(t_M), \mathcal{N}_{F,0}(|\psi(t_M)|)\psi(t_M))_{\s^2_x}|\leq (\frac{1}{M^{\frac{(n-1)\beta}{2}}}+\frac{1}{M})C(\sup\limits_{t}\|\psi(t)\|_{H^1})\leq \frac{-E}{100}\|\psi_s(x)\|_{\s^2_x}^2
\eeq
and
\eq
|(\psi(t_M), g( t_M)^{-2}V(x/g( t_M))\psi(t_M))|\lesssim \langle t_M\rangle^{-2\epsilon}\|V(x)\|_{\s^\infty_x}\|\psi(t_0)\|_{\s^2_x}^2\leq \frac{-E}{100}\|\psi_s(x)\|_{\s^2_x}^2
\eeq
which implies that for all $M\geq \tilde{M}$, $t_M\geq t_0$, $t_0$ sufficiently large,
\eq
(\psi(t_M), (H_0+g( t_M)^{-2}V(\frac{x}{g(t_M)})+2\mathcal{N}_{F,0}(|\psi(t_M)|))\psi(t_M))_{\s^2_x}\geq \frac{E}{50}\| \psi_s(x)\|_{\s^2_x}^2\label{eqrside}
\eeq
since
\eq
(\psi(t_M),H_0\psi(t_M))_{\s^2_x}\geq 0.
\eeq
Based on \eqref{mainineq}, \eqref{eq103} and \eqref{eqrside}, we obtain a contradiction, which finishes the proof.
\end{proof}

\subsection{A typical focusing example}
A concrete nonlinearity to which Theorem \ref{thm2} applies is
\eq
\mathcal{N}(k)=-\lambda \frac{k}{1+k^2}.
\eeq
\begin{lemma}\label{typicalNP}Suppose $V$ satisfies \eqref{Linear:con} and the $V$-part of \eqref{W} (i.e., $x\cdot\nabla V\in \s^\infty_x$); suppose $H=H_0+V$ satisfies \eqref{eq: H lambda} and $g$ satisfies \eqref{g(t)}. If
\eq
\mathcal{N}(k)=-\lambda \frac{k}{1+k^2}
\eeq
for $\lambda>0$ sufficiently large, then the hypotheses of Theorem \ref{thm2} are satisfied, and the corresponding solution exhibits at least two bubbles as $t\to\infty$.
\end{lemma}
\begin{proof}For $\mathcal N(k)=-\lambda k/(1+k^2)$, a direct computation gives $\mathcal N_{F,0}(k)=\int_0^k q\mathcal N(q)/k^2\,dq=-\lambda(k-\arctan k)/k^2$ (for $k>0$). By the standard concentration-compactness/variational argument for saturated focusing nonlinearities, for $\lambda>0$ sufficiently large the stationary equation $(H_0+2\mathcal N_{F,0}(|\psi_s|))\psi_s=E\psi_s$ admits a radial $H^1$ solution with $E<0$; equivalently, there is a soliton $\psi_s(x)$ to \eqref{NF00} satisfying \eqref{soliton}.

Taking $\beta=1$, the bound $|\mathcal N(k)|=\lambda k/(1+k^2)\leq \lambda k$ verifies \eqref{NFandN}. The additional regularity used by Lemma~\ref{lem:NPH1} is also direct: the saturation $|\mathcal N(k)|\leq \lambda$ gives $\mathcal N\in \s^\infty_x(\mathbb R_+)$, and $\mathcal N'(q)=-\lambda(1-q^2)/(1+q^2)^2$ yields
\eq
|\mathcal N_F(k)|=\Bigl|\int_0^k q^2\mathcal N'(q)\,dq\Bigr|\lesssim \lambda\!\int_0^k\!\frac{q^2}{1+q^2}\,dq\lesssim k^2.
\eeq
Together with the standing hypotheses on $V$, $H$, $g$, Lemma~\ref{lem:NPH1} therefore applies and produces a global $H^1$ solution $\psi\in C([t_0,\infty);H^1(\mathbb R^n))$ of \eqref{NP} satisfying \eqref{NP:H1bound}.

Thus all hypotheses of Theorem~\ref{thm2} are met, and the solution $\psi(t)$ to
\eq
\begin{cases}
i\partial_t\psi(t)=(H_0+g(t)^{-2}V(x/g(t))-\lambda \frac{|\psi(t)|}{1+|\psi(t)|^2})\psi(t)\\
\psi(x,t_0)=\psi_s(x)+e^{-iD\ln(g(t_0))}\psi_b(x)
\end{cases},
\eeq
exhibits at least two bubbles with different patterns. We finish the proof.
\end{proof}

\bibliographystyle{uncrt}

\end{document}